\documentclass[preprint, 5p, twocolumn]{elsarticle}
\journal{Journal of Terramechanics}

\usepackage{amsmath}
\usepackage{amssymb}
\usepackage{amsthm}

\usepackage{booktabs}
\usepackage{multirow}

\usepackage{tikz}
\usepackage{graphicx}
\usepackage{subcaption}

\usepackage{algorithm}
\usepackage{algpseudocode}

\usetikzlibrary{shapes.geometric}

\newcommand{\norm}[1]{\lVert#1\rVert}

\begin{document}

\begin{frontmatter}

\title{Tensor-Train Compression of Discrete Element Method Simulation Data\tnoteref{t1}}
\tnotetext[t1]{DISTRIBUTION A. Approved for public release; distribution unlimited. OPSEC $\# 6887$.}

\author[1]{Saibal De}
\ead{sde@sandia.gov}

\author[2]{Eduardo Corona}
\ead{eduardo.corona@colorado.edu}

\author[3]{Paramsothy Jayakumar}
\ead{paramsothy.jayakumar.civ@army.mil}

\author[4]{Shravan Veerapaneni}
\ead{shravan@umich.edu}

\address[1]{University of Michigan, Ann Arbor; now at Sandia National Laboratories, Livermore}
\address[2]{University of Colorado, Boulder}
\address[3]{US Army DEVCOM Ground Vehicle Systems Center (GVSC)}
\address[4]{University of Michigan, Ann Arbor}

\begin{abstract}
    We propose a framework for discrete scientific data compression based on the tensor-train (TT) decomposition. Our approach is tailored to handle unstructured output data from discrete element method (DEM) simulations, demonstrating its effectiveness in compressing both raw (e.g.\ particle position and velocity) and derived (e.g.\ stress and strain) datasets. We show that geometry-driven ``tensorization'' coupled with the TT decomposition (known as quantized TT) yields a hierarchical compression scheme, achieving high compression ratios for key variables in these DEM datasets.
\end{abstract}

\end{frontmatter}

\section{Introduction}

\begin{figure*}
    \centering
    
    \includegraphics[width=0.23\textwidth]{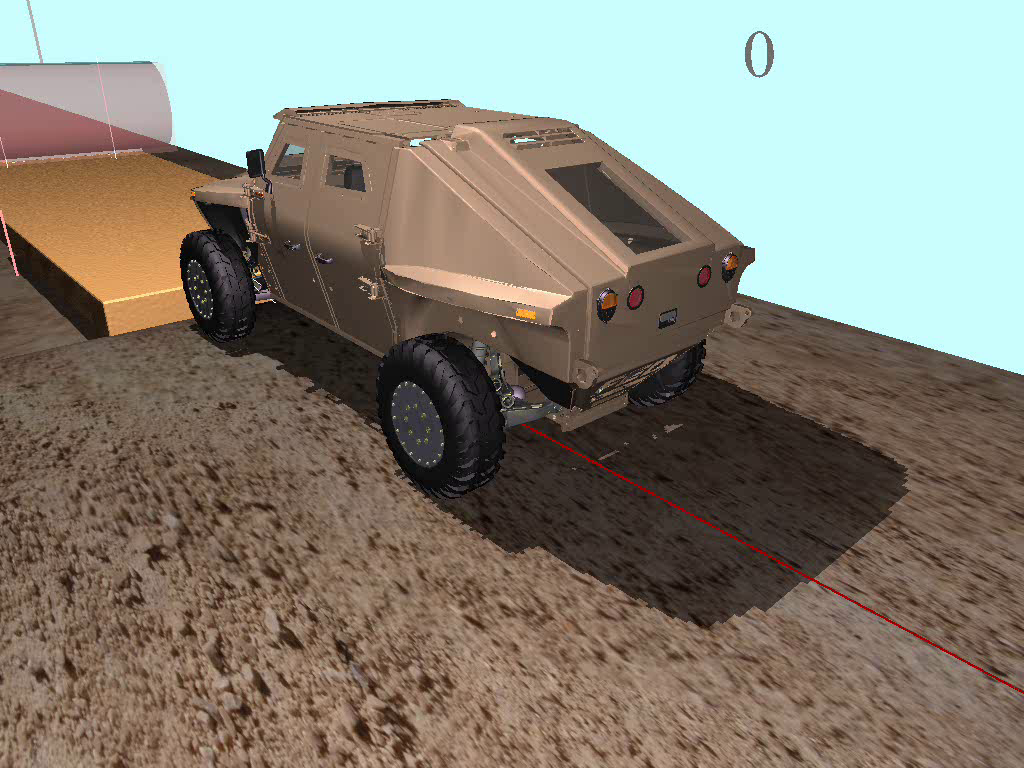}
    ~
    \includegraphics[width=0.23\textwidth]{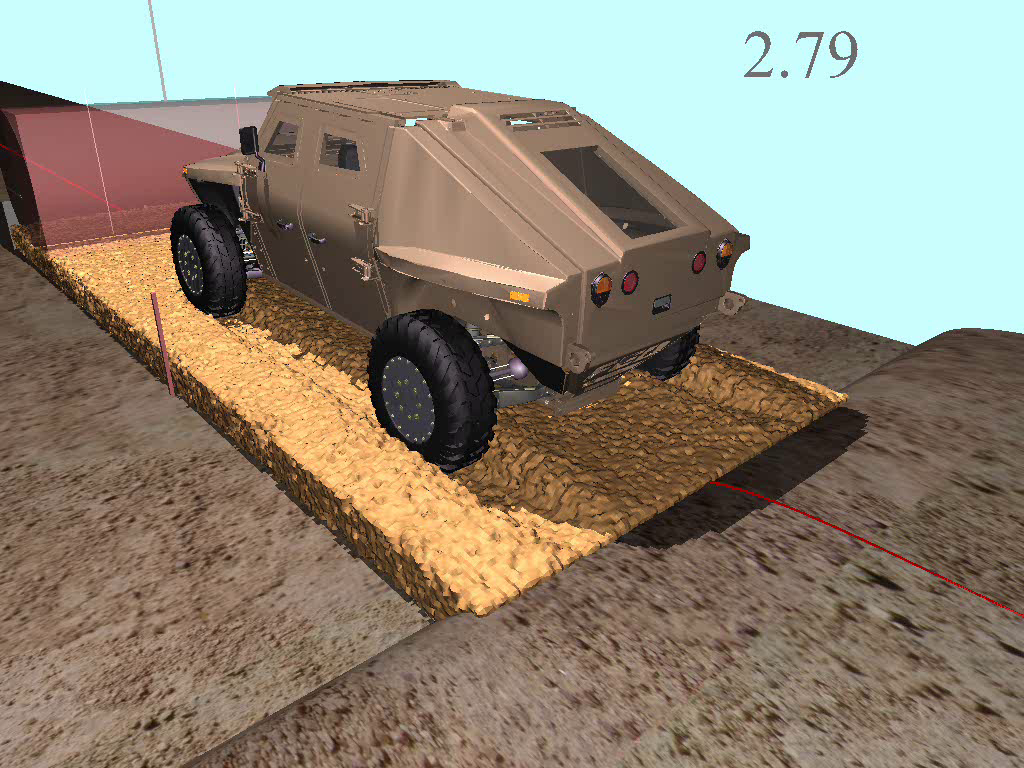}
    ~
    \includegraphics[width=0.23\textwidth]{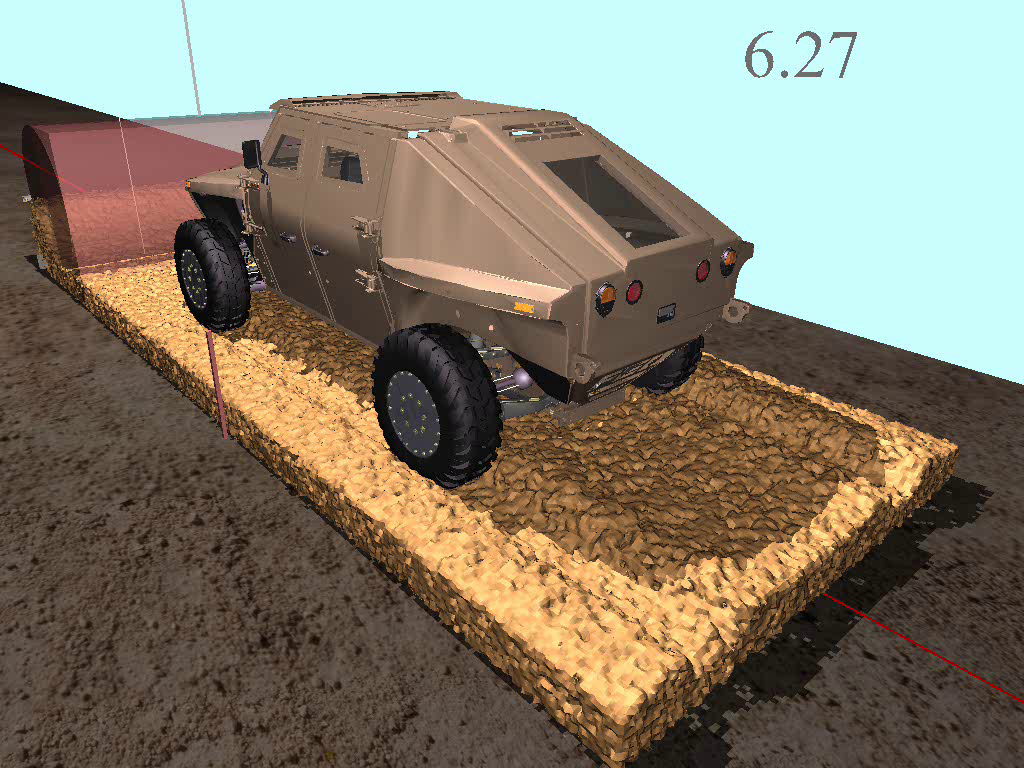}
    ~
    \includegraphics[width=0.23\textwidth]{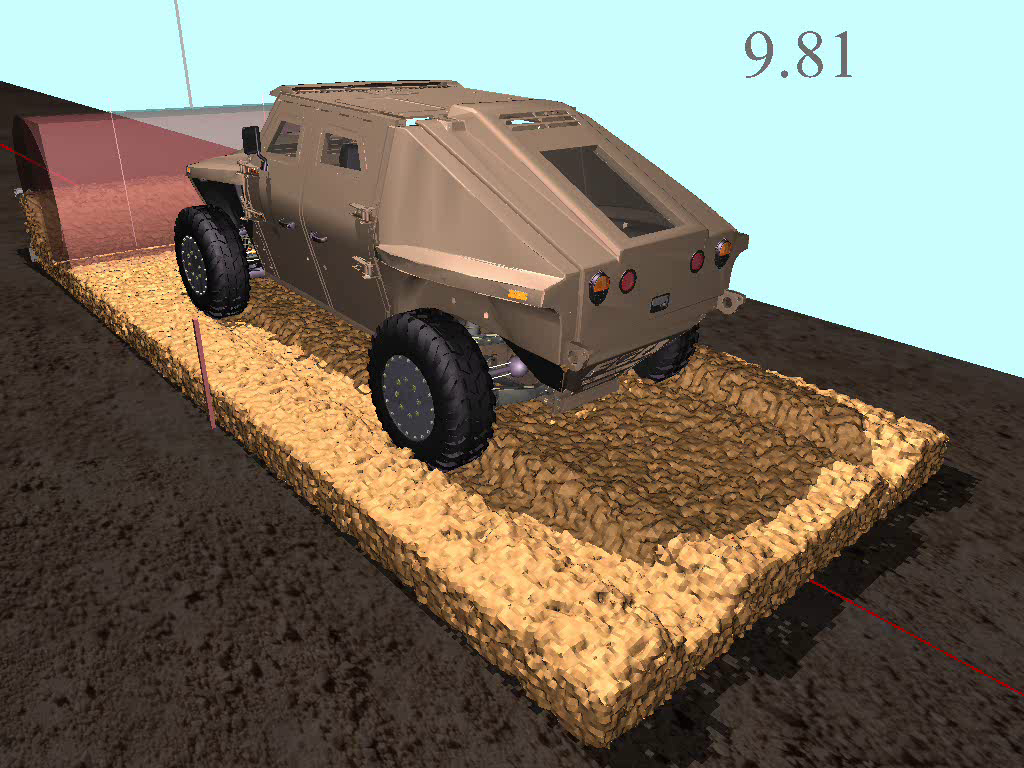}
    
    \caption{{\em Snapshots from a DEM simulation of a vehicle traversing an off-road terrain, with the soil patch modeled as a collection 488,376 spherical particles. The position and velocity of DEM soil particles is stored at each of the 1,533 timesteps with an uniform step size, $\Delta t = 0.03$~s. The resulting datasets are approximately 16.7~GB each for position and velocity, when recorded in IEEE 64-bit floating point representation. The techniques developed in this work compress the position data by a factor of $1.8 \times 10^7$ to less than 1~KB, and the velocity data by a factor of $3.1 \times 10^3$ to approximately 5.2~MB, while guaranteeing that normalized root mean squared error (nRMSE) is smaller than $10^{-1}$ (see Section \ref{sec:soilvehicle} for more details). }}
    \label{fig:soilvehicle-setup}
\end{figure*}

The discrete element method (DEM) is widely recognized as an effective simulation framework for granular media modeled as a collection of individual particles. It uses simple first-principle physical laws to describe the particle response to body forces (e.g.\ gravity and inertia) and dissipative contact forces (e.g.\ friction), leading to high-fidelity and robust predictions. In recent years, substantial progress has been made to address the algorithmic challenges in DEM simulations: in fast optimization solvers \cite{stewart2000rigid, tasora2010DEMC, melanz2017comparison, corona2018tensor}, parallel computing frameworks \cite{coumans2013bullet, mazhar2013chrono, negrut2014parallel, berger2015hybrid, de2019scalable, corona2020scalable}, and multiscale modeling \cite{li2016advances, zhao2017hierarchical, sugiyama2017high, yamashita2020parallelized}. Unfortunately, for large-scale problems such as those in vehicle-terrain simulation (see Figure~\ref{fig:soilvehicle-setup}), each run can still take days or even weeks to complete. In learning scenarios such as mobility map construction, each of these expensive runs then corresponds to only one point in the exploration of a high-dimensional parametric space \cite{jayakumar2019efficient, marple2020active}.

As algorithmic and learning frameworks continue to improve, addressing large dataset storage and management problems is imperative. DEM vehicle-soil simulations, for example, generate a lot of data; this might include particle positions and velocities, pairwise contact forces and soil-vehicle interaction forces at each timestep. One simulation can easily produce tens of gigabytes of data. In current state-of-the-art implementations of learning tasks, only the bare minimum of the computational results is kept. This is not only wasteful of weeks or months of expensive, energy-intensive computational work, but it also precludes any current or future knowledge we might obtain from the entire dataset.

We note that these issues are far from exclusive to DEM. Across various disciplines that employ physics-based simulations, tremendous amounts of data is produced, often outpacing the current storage capacities. This limits our ability to learn, and makes data transfer, visualization and analysis significant bottlenecks in our design pipelines. This is most acutely felt in online, interactive environments in which streaming data must be processed under restrictive computational and storage capabilities. Data compression techniques are therefore of great interest across the entire computational simulation landscape.

\subsection{Related Works}

There is a vast literature on general scientific data reduction methods. While general purpose encoders such as \texttt{fpzip} \cite{lindstrom2006fast} have achieved a degree of success exploiting data coherence in scientific datasets, more efficient representations require us to incorporate assumptions on specific types of data-sparsity or inherent physical structures we want to target. Techniques such as principal component analysis (PCA), wavelet transforms and multi-resolution schemes have been traditionally deployed to exploit inherent structure due to smoothness and data locality. For a more complete and informative survey on general scientific data compression, see \cite{li2018data}.

Scientific data, such as those obtained in DEM simulations, is almost inevitably multi-dimensional, as variables of interest are tracked in space, time, and across a range of parameters. Tensor factorization schemes are a natural extension of matrix-factorization based data compression methods like the PCA. They explicitly tackle the high-dimensional nature of scientific data. Key features of tensor decompositions include addressing the curse of dimensionality and producing significantly reduced representations of datasets that allow rapid full or partial reconstruction and fast linear algebra in compressed format at the expense of allowing user-specified error tolerance.

Due to recent algorithmic developments, there is a rapidly growing body of literature for tensor decomposition algorithms and their applications to data mining, data fusion and scientific data compression \cite{kolda2009review,acar2008unsupervised,grasedyck2013literature,cichocki2016tensor}. The work reviewed in these articles suggests the various tensor representations (canonical polyadic, Tucker, tensor-train) and their hierarchical counterparts ($\mathcal{H}$-Tucker, quantized tensor-train) can be incredibly effective as tools for data compression and analysis. For instance, in \cite{austin2016parallel}, the authors developed a distributed implementation of the Tucker decomposition, demonstrating dramatic compression of combustion datasets (e.g. 760-200K fold compression for a $10^{-2}$ relative tolerance) and fast partial reconstruction in limited computing environments. More recently, there has been a large push towards developing tensor compression techniques for streaming data with canonical polyadic \cite{nguyen2016fast, smith2018streaming}, Tucker \cite{malik2018low, sun2020low} and tensor-train \cite{liu2018incremental, liu2020scalable, abed2021adaptive} factorizations.

In the context of DEM simulations, methods exploiting low-rank data structure have not focused on data compression and reconstruction, but on producing reduced order models (ROM). Many of these approaches still leverage low-rank factorizations, such as the singular value decomposition (SVD), to construct the reduced representation. For instance, \cite{zhong2018adaptive} uses the related proper orthogonal decomposition (POD) method to extract the principal modes of the dynamics from time snapshots of a damped DEM simulation.

The reliance of ROM methods on low rank matrix techniques to process DEM training datasets, and their success at yielding meaningful bases for the dynamics, is strong evidence for the potential of matrix and tensor compression methods like the one proposed in this work.

\subsection{Our Contributions}

We tailor the tensor-train (TT) format \cite{oseledets2011tensor} to compress both raw and derived data from DEM simulations. Our contributions to the state-of-the-art of scientific data compression are as follows:

\begin{itemize}
    \item We present what is, to the best of our knowledge, the first adaptation of tensor compression methods to nonsmooth granular media simulation data.

    \item Our methods feature a compression scheme for streaming data, in which batches of simulation data can be efficiently incorporated to a compressed representation as soon as they are generated.

    \item Most applications of tensor compression methods in the literature focus on compression along the dimensions suggested by the data (e.g. spatial $x$, $y$, and $z$ directions, time). By ``tensorizing'' input data using a standard space-filling curve ordering, we demonstrate the ability of quantized tensor-train (QTT) methods to provide further compression by exploiting redundancies in the data at different levels of resolution.
\end{itemize}

We demonstrate high compression ratios (1K-10M fold reductions) for raw and derived DEM datasets while keeping the relative reconstruction error smaller than $10^{-1}$. This enables us to store several gigabytes of simulation data in few megabytes or kilobytes of disk space, a small fraction of the original dataset size. Fast algorithms for TT array access and arithmetic \cite{oseledets2011tensor} may then be used to perform post-processing and learning tasks in compressed TT form.

\section{The Tensor-Train Decomposition}

From a computational perspective, tensors are a direct generalization of matrices in higher than two dimensions: they are multi-dimensional arrays of real numbers. The entries of a $d$-dimensional tensor $\mathcal{X} \in \mathbb{R}^{n_1 \times \cdots \times n_d}$ are expressed as $\mathcal{X}(i_1, \ldots, i_d)$ with indices $1 \leq i_k \leq n_k$ for $1 \leq k \leq d$. Given $n = \max\{n_1, \ldots, n_d\}$ and assuming tensor dimensions are comparable in size, it can be readily observed that storing all of its entries incurs $O(n^d)$ storage cost. Additionally, any operation accessing all of the entries of the tensor will have a computational complexity that is similarly exponential in the dimensionality $d$.

Low-rank tensor decompositions leverage multilinear structure of the underlying tensor --- in essence generalizing the idea of low-rank matrix factorizations --- in order to overcome this ``curse of dimensionality''. Several formats for such tensor decompositions have been developed: in \cite{kolda2009review}, the authors provide a detailed review of the canonical polyadic (CP) factorization and the Tucker decomposition, and \cite{oseledets2011tensor} introduces the tensor-train (TT) format.

In this work, we use the TT decomposition as our main tool for data compression. This format does not have a high-dimensional core tensor (unlike the Tucker factorization), and it can be computed in a numerically stable manner with provable optimality bounds (unlike the CP decomposition) via a sequence of low-rank matrix factorizations. Additionally, the quantized tensor-train (QTT) provides a straight-forward method for hierarchical multiresolution compression of datasets; a generalization of this idea is the basis for the $\mathcal{H}$-Tucker decomposition.

In this section, we briefly review the TT format; see \cite{oseledets2011tensor} for a more detailed exposition. The TT decomposition of a $d$-dimensional tensor $\mathcal{X} \in \mathbb{R}^{n_1 \times \cdots \times n_d}$ is given by
\begin{equation}
    \label{eq:tt-format}
    \begin{split}
        &\mathcal{X}(i_1, \ldots, i_d) = \\
        &\qquad \sum_{\alpha_0 = 1}^{r_0} \cdots \sum_{\alpha_d = 1}^{r_d} \mathcal{X}_1(\alpha_0, i_1, \alpha_1) \cdots \mathcal{X}_d(\alpha_{d - 1}, i_d, \alpha_d).
    \end{split}
\end{equation}
The limits $r_0, \ldots, r_d$ of the sums are known as the TT ranks of the tensor --- they are a measure of the complexity of the underlying data. The three-dimensional factors $\mathcal{X}_k \in \mathbb{R}^{r_{k - 1} \times n_k \times r_k}$ for $1 \leq k \leq d$ are known as TT cores. The boundary ranks are generally assumed to be singletons ($r_0 = r_d = 1$); in this setup, a useful reformulation of \eqref{eq:tt-format} is the matrix product state
\begin{equation}
    \mathcal{X}(i_1, \ldots, i_d) = X_1(i_1) \cdots X_d(i_d),
\end{equation}
where the matrices on the right-hand side are slices of the TT cores:
\begin{equation}
    X_k(i_k) = \mathcal{X}_k(:, i_k, :) \in \mathbb{R}^{r_{k - 1} \times r_k}.
\end{equation}

The advantage of using the TT format in terms of storage complexity is immediate: given a tensor with maximum TT rank $r = \max\{r_1, \ldots, r_{d - 1}\}$, storing the TT cores requires $O(d n r^2)$ memory. A large array of tensors relevant to practical applications can be well approximated by tensors with $r \ll \sqrt{n^d}$, which leads to significant storage savings \cite{oseledets2011tensor2}.

The TT decomposition of the tensor $\mathcal{X}$ in effect corresponds to providing low-rank matrix factorizations for the so-called unfolding matrices: the $k$-th unfolding of the $d$-dimensional tensor $\mathcal{X}$ constructs a matrix by grouping the first $k$ indices of the tensor into rows, and the remaining indices into columns
\begin{equation}
    \operatorname{unfold}_k(\mathcal{X})(\overline{i_1, \ldots, i_k}, \overline{i_{k + 1}, \ldots, i_d}) = \mathcal{X}(i_1, \ldots, i_d)
\end{equation}
where
\begin{align}
    \label{eq:long_index}
    \begin{split}
        \overline{i_1, \ldots, i_k}       &= 1 + \sum_{j = 1}^k \left(\prod_{j' = 1}^{j - 1} n_{j'}\right) (i_j - 1), \\
        \overline{i_{k + 1}, \ldots, i_d} &= 1 + \sum_{j = k + 1}^d \left(\prod_{j' = k + 1}^{j - 1} n_{j'}\right) (i_j - 1),
    \end{split}
\end{align}
are long indices. If the tensor $\mathcal{X}$ is stored in column-major (Fortran) order, this is a reshaping of the array:
\begin{equation}
    \operatorname{unfold}_k(\mathcal{X}) = \mathsf{reshape}(\mathcal{X}, [n_1 \cdots n_k, n_{k + 1} \cdots n_d])
\end{equation}
for $1 \leq k \leq d - 1$.

Note that a low-rank factorization of $\operatorname{unfold}_k(\mathcal{X}) \approx U V$ produces a sum of products with factors depending on only one of these two long indices, while introducing an auxiliary parameter ranging from $1$ to $r_k$. We may then apply this same procedure to the factors $U$ and $V$, iterating until we get a decomposition where each factor depends on only one tensor dimension. This suggests computing TT cores using a sequence of low rank factorizations of these unfolding matrices; we obtain cores $\mathcal{Y}_1, \ldots, \mathcal{Y}_d$ of a tensor $\mathcal{Y}$ approximating the original tensor $\mathcal{X}$ such that
\begin{equation}
    \label{eq:approx_TT}
    \| \mathcal{X} - \mathcal{Y} \|_F \leq \tau_\text{relFrob} \| \mathcal{X} \|_F
\end{equation}
for some given relative Frobenius error tolerance $\tau_\text{relFrob}$. A near-optimal factorization may be computed by applying a sequence of truncated SVDs. The details of this computation are presented in Algorithm~\ref{alg:tt-svd}.

\begin{algorithm}[t]
    \caption{TT-SVD \cite{oseledets2011tensor}}
    \label{alg:tt-svd}
    \begin{algorithmic}[1]
        \Require {Tensor $\mathcal{X} \in \mathbb{R}^{n_1 \times \cdots \times n_d}$, relative tolerance $\tau_\text{relFrob}$}
        \Ensure {Cores $\mathcal{Y}_1, \ldots, \mathcal{Y}_d$ s.t.\ $\norm{\mathcal{X} - \mathcal{Y}}_F \leq \tau_\text{relFrob} \norm{\mathcal{X}}_F$}
        \Statex
        \State {Compute truncation: $\delta \gets \tau_\text{relFrob} \norm{\mathcal{X}}_F / \sqrt{d - 1}$}
        \State {Initialize $M \gets \textsf{reshape}(\mathcal{X}, [n_1, n_2 \cdots n_d])$}
        \State {Initialize $r_0 \gets 1$}
        \State {Initialize $r_d \gets 1$}
        \For {$k = 1, \ldots, d - 1$}
            \State {$U, \Sigma, V \gets \textsf{truncated\_svd}_\delta(M)$}
            \State {$r_k \gets \textsf{rank}(\Sigma)$}
            \State {$\mathcal{Y}_k \gets \textsf{reshape}(U, [r_{k - 1}, n_k, r_k])$}
            \State {$M \gets \textsf{reshape}(\Sigma V^\top, [r_k n_{k + 1}, n_{k + 2} \cdots n_d])$}
        \EndFor
        \State {$\mathcal{Y}_d \gets \textsf{reshape}(M, [r_{d - 1}, n_d, r_d])$}
    \end{algorithmic}
\end{algorithm}

The TT format is particularly suitable for linear algebra \cite{oseledets2011tensor}; in particular, we can add and subtract two tensors in TT format without needing to compute their full representations. These fast arithmetic operations usually inflate the TT ranks, leading to sub-optimal representations. Fortunately, we can recompress these tensors by applying matrix factorizations to each of the TT cores \cite{oseledets2011tensor}. This TT-rounding process constructs from input TT tensor $\mathcal{X}$ and relative Frobenius error tolerance $\tau_\text{relFrob}$ a new TT tensor $\mathcal{Y}$ satisfying \eqref{eq:approx_TT}. The steps for TT-rounding are included in Algorithm~\ref{alg:tt-round}; we note that the Frobenius norm of the input tensor $\mathcal{X}$ in TT-format can be computed directly from the tensor cores, and we do not need to reconstruct the full tensor \cite{oseledets2011tensor}. This process has computational complexity of $O(d n r^3)$, where $r$ is the maximal TT rank of the sub-optimal tensor.

\begin{algorithm}[t]
    \caption{TT-rounding \cite{oseledets2011tensor}}
    \label{alg:tt-round}
    \begin{algorithmic}[1]
        \Require {TT cores $\mathcal{X}_1, \ldots, \mathcal{X}_d$, relative tolerance $\tau_\text{relFrob}$}
        \Ensure {Cores $\mathcal{Y}_1, \ldots, \mathcal{Y}_d$ s.t.\  $\norm{\mathcal{X} - \mathcal{Y}}_F \leq \tau_\text{relFrob} \norm{\mathcal{X}}_F$}
        \Statex
        \State {Compute truncation: $\delta \gets \tau_\text{relFrob} \norm{\mathcal{X}}_F / \sqrt{d - 1}$}
        \State {Initialize cores $\mathcal{Y}_k \gets \mathcal{X}_k$ for $1 \leq k \leq d$}
        \For {$k = d, d - 1, \ldots, 2$}
            \State {$M \gets \mathsf{reshape}(\mathcal{Y}_k, [r_{k - 1}, n_k r_k])$}
            \State {$Q, R \gets \mathsf{qr}(M^\top)$}
            \State {$r_{k - 1}' \gets \mathsf{rank}(R)$}
            \State {$\mathcal{Y}_k \gets \mathsf{reshape}(Q^\top, [r_{k - 1}', n_k, r_k])$}
            \State {$M \gets \mathsf{reshape}(\mathcal{Y}_{k - 1}, [r_{k - 2} n_{k - 1}, r_{k - 1}])$}
            \State {$\mathcal{Y}_{k - 1} \gets \mathsf{reshape}(M R^\top, [r_{k - 2}, n_{k - 1}, r_{k - 1}'])$}
            \State {$r_{k - 1} \gets r_{k - 1}'$}
        \EndFor
        \For {$k = 1, \ldots, d - 1$}
            \State {$M \gets \mathsf{reshape}(\mathcal{Y}_k, [r_{k - 1} n_k, r_k])$}
            \State {$U, \Sigma, V \gets \mathsf{svd}_\delta(M)$}
            \State {$r_k' \gets \mathsf{rank}(\Sigma)$}
            \State {$\mathcal{Y}_k \gets \mathsf{reshape}(U, [r_{k - 1}, n_k, r_k'])$}
            \State {$M \gets \mathsf{reshape}(\mathcal{Y}_{k + 1}, [r_k, n_{k + 1} r_{k + 1}])$}
            \State {$\mathcal{Y}_{k + 1} \gets \mathsf{reshape}(\Sigma V^\top M, [r_k', n_{k + 1}, r_{k + 1}])$}
            \State {$r_k \gets r_k'$}
        \EndFor
    \end{algorithmic}
\end{algorithm}

\section{Tensor-Train for DEM Data Compression}

In this section, we propose an adaptation of the TT algorithms outlined in Section 2 to compress scientific datasets. For this purpose, we review a QTT approach for hierarchical domain decomposition of \emph{tensorized} vectors and matrices obtained from function sampling. This method reliably produces low-rank TT factorizations when reordering the data in a spatially coherent fashion. Using this observation, we then extend the QTT method to tackle the specific challenges present in data generated by DEM simulations.

\subsection{Tensorization}

One of the main advantages of the TT decomposition is the linear scaling of storage complexity $\mathcal{O}(d n r^2)$ with the dimensionality $d$ of the data tensor. This makes it ideal for compressing very high-dimensional tensors. We take full advantage of this aspect of TT factorization by systematically increasing the dimensionality of the data (a.k.a.\ tensorizing the data) before applying the compression algorithms.

This approach is motivated by the hierarchical domain decomposition techniques used in partial differential equation solvers, especially in the context of boundary integral equation (BIE) methods \cite{greengard1987fast, ho2014fast, ho2015hierarchical_ie, ho2015hierarchical_de, corona2017tensorbie}. In the following, we first present two examples which demonstrate the advantage of tensorization in the context of structured data --- where we sample functions to construct the datasets. We note these functions are smooth, but develop a localized discontinuity as the function parameter goes to zero. We then extend this approach to the case of unstructured particle data from DEM simulations.

\subsubsection{Compression of a Univariate Function}

We consider the function
\begin{equation}
    f_\delta(x) = \ln\frac{1}{| x - 1/2 | + \delta}, \quad \delta > 0,
\end{equation}
which is smooth everywhere, but its peak at $x = 1/2$ gets sharper as we decrease $\delta \to 0$. We sample this function uniformly in the interval $(0, 1)$ for a fixed value of $\delta$ and construct a vector $F_\delta$ of length $2^d$ and entries
\begin{equation}
    F_\delta(i) = f_\delta\left(\frac{i - 1/2}{2^d}\right), \quad 1 \leq i \leq 2^d.
\end{equation}
We tensorize this vector by constructing a binary tree of depth $d$ and assigning the entries of the vector $F_\delta$ to the leaves of the binary tree. This creates a $d$-dimensional tensor of size\footnote{For brevity, we adapt the exponentiation notation to describe tensor sizes; thus the size of a $2 \times 2 \times 2 \times 2$ tensor is shortened to $2^{\times 4}$, and that of a $4 \times 3 \times 3 \times 4 \times 4$ tensor is shortened to $4 \times 3^{\times 2} \times 4^{\times 2}$.} $2^{\times d}$ whose entries are given by
\begin{equation}
    \mathcal{F}_\delta(i_1, \ldots, i_d) = F_\delta(\overline{i_1, \ldots, i_d}), \quad 1 \leq i_1, \ldots, i_d \leq 2
\end{equation}
where $\overline{i_1, \ldots, i_d}$ is a long index as defined in \eqref{eq:long_index}.

\begin{figure}[t]
    \centering
    \begin{tikzpicture}[scale=0.8]
    \draw (5, 3) -- node [anchor=east] {$i_4 = 1\ \ $} (3, 2);
    \draw (5, 3) -- node [anchor=west] {$\ \ 2$}(7, 2);

    \foreach \x in {3, 7} {
        \draw ({\x}, 2) -- node [anchor=east] {$\ifthenelse{\x = 3}{i_3 = }{}1\ $} (\x - 1, 1);
        \draw ({\x}, 2) -- node [anchor=west] {$\ 2$} (\x + 1, 1);
    }

    \foreach \x in {2, 4, 6, 8} {
        \draw (\x - 1, -1.5) -- (\x + 1, -1.5) -- (\x + 1, 1) -- (\x - 1, 1);
    }
    \draw (1, -1.5) -- (1, 1);

    \foreach \x in {2, 4, 6, 8} {
        \foreach \y in {-2, -1, 0, 1} {
            \draw[draw=gray, fill=gray] (\x, 0.5 * \y) circle (0.09);
        }
        \draw [thick, gray] (\x, 0.5) -- (\x, -1.0);
    }
    \foreach \x in {2, 4, 6} {
        \draw [thick, gray] (\x, -1.0) -- (\x + 2, 0.5);
    }

    \foreach \x in {1, 5, ..., 16} {
        \node [anchor=west] at (0.5 * \x + 1.5,  0.5) {\scriptsize $F(\x)$};
    }
    \node [anchor=east] at (1, 0.5) {$\overline{i_1, i_2} = 1$};

    \foreach \x in {2, 6, ..., 16} {
        \node [anchor=west] at (0.5 * \x + 1,  0.0) {\scriptsize $F(\x)$};
    }
    \node [anchor=east] at (1, 0.0) {$\overline{i_1, i_2} = 2$};

    \foreach \x in {3, 7, ..., 16} {
        \node [anchor=east] at (0.5 * \x + 0.5, -0.5) {\scriptsize $F(\x)$};
    }
    \node [anchor=east] at (1, -0.5) {$\overline{i_1, i_2} = 3$};

    \foreach \x in {4, 8, ..., 16} {
        \node [anchor=east] at (0.5 * \x, -1.0) {\scriptsize $F(\x)$};
    }
    \node [anchor=east] at (1, -1.0) {$\overline{i_1, i_2} = 4$};
\end{tikzpicture}
    \caption{{\em Tensorization of a vector $F(i)$ with 16 elements into a three-dimensional tensor $\mathcal{F}^{\{3\}}(\overline{i_1, i_2}, i_3, i_4)$ of size $4 \times 2 \times 2$ via a truncated binary tree construction. The tensor indices $i_3$ and $i_4$ represent whether we go left or right when moving down the levels of the tree, and the long index $\overline{i_1, i_2}$ identifies the vector entry within the four-element leaf. For instance, the entry $F(11)$ is represented as $\mathcal{F}^{\{3\}}(3, 1, 2)$. Assuming a column-major order for multidimensional arrays, this is essentially a reshaping: $\mathcal{F}^{\{3\}} = \mathsf{reshape}(F, [4, 2, 2])$.}}
    \label{fig:tensorization}

    \includegraphics[width=\linewidth]{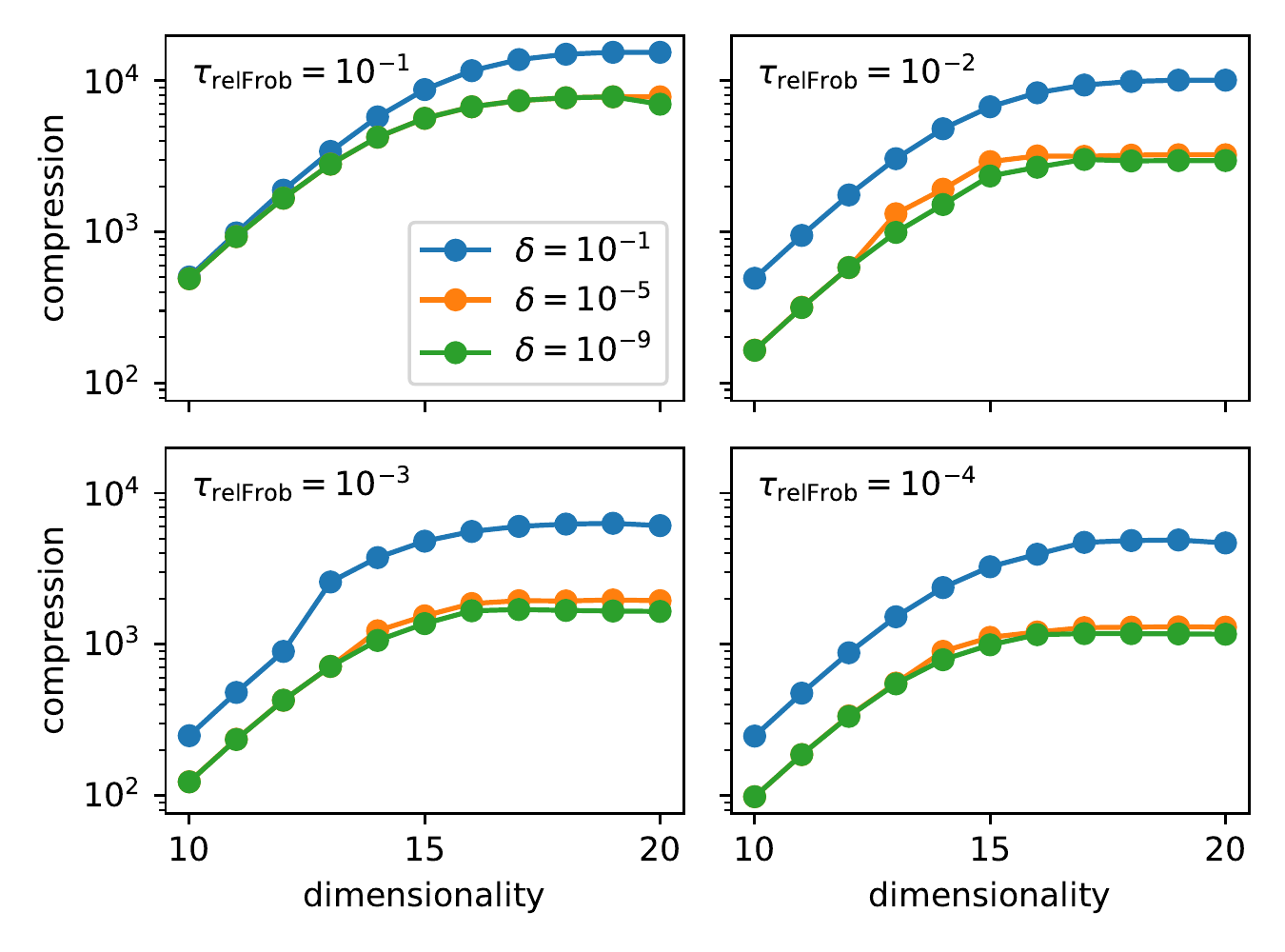}
    \caption{{\em Compression ratios of the tensorized function samples $\mathcal{F}_\delta^{\{l\}}$ for different values of smoothness parameter $\delta$ and tensorization level $l$ of the reshaped dataset. As we increase $l$, TT compression ratio steadily improves and eventually flattens out. As we decrease $\delta$, making the underlying function more sharper near $x = 1/2$, the compression ratio decreases. These observations hold true irrespective of the value of the relative Frobenius error tolerance parameter $\tau_\text{relFrob}$.}}
    \label{fig:demo_function}
\end{figure}

More generally, for a fixed value of $1 \leq l \leq d$, we can stop at the $(l - 1)$-th level of the binary tree and collect $2^{d - l + 1}$ entries of the vector per leaf. This creates a $l$-dimensional view of the vector
\begin{equation}
    \mathcal{F}_\delta^{\{l\}}(\overline{i_1, \ldots, i_{d - l + 1}}, i_{d - l + 2}, \ldots, i_d) = F_\delta(\overline{i_1,\cdots,i_d})
\end{equation}
where the indices $i_k$ corresponds to which branch to follow as we move down the levels of the binary tree from the root to the appropriate leaf. We demonstrate this process in Figure~\ref{fig:tensorization} for $d = 4$ and $l = 3$.

When compressing the tensors $\mathcal{F}_\delta^{\{l\}}$ for different values of $l$, the first low-rank decomposition in the TT-SVD algorithm is performed on different unfolding matrices
\begin{equation}
    \operatorname{unfold}_1(\mathcal{F}_\delta^{\{l\}}) = \textsf{reshape}(\mathcal{F}_\delta, [2^{d - l + 1}, 2^{l - 1}])
\end{equation}
of the fully resolved tensor $\mathcal{F}_\delta \equiv \mathcal{F}_\delta^{\{d\}} \in \mathbb{R}^{2^{\times d}}$. Intuitively, increasing $l$ leads to decreasing the number of rows in the first unfolding matrices. This in effect allows us to leverage the locally smooth behavior of the underlying function in most of the function domain more effectively.

We apply this tensorization approach to the function samples for $d = 20$ and for $l = 10, \ldots, 20$. We compress the resulting $l$-dimensional tensors $\mathcal{F}_\delta^{\{l\}}$ for different values of the sharpness parameter $\delta \in \{10^{-1}, 10^{-5}, 10^{-9}\}$ and the relative Frobenius error tolerance $\tau_\text{relFrob} \in \{10^{-1}, 10^{-2}, 10^{-3}, 10^{-4}\}$ of the TT-SVD algorithm. We plot the resulting compression ratios against the number of levels $l$ in Figure~\ref{fig:demo_function}. We note irrespective of the values of $\delta$ and $\tau_\text{relFrob}$, the compression ratio increases as we increase the dimensionality of the tensor.

We also note that while all vectors are highly compressible, as we decrease $\delta$, making the function more sharply peaked, the compression ratio decreases. This suggests that the presence of localized discontinuities in a dataset may limit the compression ratio achievable at a given target accuracy.

\subsubsection{Compression of a Bivariate Function}

\begin{figure}[t]
    \centering
    \begin{tikzpicture}[scale=0.4]
    \foreach \i in {0, 1, ..., 15} {
        \foreach \j in {0, 1, ..., 15} {
            \draw [draw=gray, fill=gray] (\i, \j) circle (0.15);
        }
    }
    
    \foreach \i in {0, 4, 8, 12} {
        \foreach \k in {0, 1, 2, 3} {
            \draw (\i + 1.5, 17) -- (\i + \k, 16);
            \draw (-2, \i + 1.5) -- (-1, \i + \k);
        }
    }
    \foreach \i in {0, 8} {
        \foreach \k in {0, 4} {
            \draw (\i + 3.5, 18) -- (\i + \k + 1.5, 17);
            \draw (-3, \i + 3.5) -- (-2, \i + \k + 1.5);
        }
    }
    \foreach \i in {0} {
        \foreach \k in {0, 8} {
            \draw (\i + 7.5, 19) -- (\i + \k + 3.5, 18);
            \draw (-4, \i + 7.5) -- (-3, \i + \k + 3.5);
        }
    }

    \node[anchor=east] at (-2, 13.5) {$\overline{i_1, i_2}$};
    \node[anchor=east] at (-3, 11.5) {$i_3$};
    \node[anchor=east] at (-4,  7.5) {$i_4$};

    \node[anchor=south east] at (1.5, 17) {$\overline{j_1, j_2}$};
    \node[anchor=south] at (3.5, 18) {$j_3$};
    \node[anchor=south] at (7.5, 19) {$j_4$};

    \draw [very thick, dashed] (7.5, -0.5) -- (7.5, 15.5);
    \draw [very thick, dashed] (-0.5, 7.5) -- (15.5, 7.5);

    \foreach \i in {0, 8} {
        \foreach \j in {0, 8} {
            \draw [very thick, dotted] (\i + 3.5, \j - 0.25) -- (\i + 3.5, \j + 7.25);
            \draw [very thick, dotted] (\i - 0.25, \j + 3.5) -- (\i + 7.25, \j + 3.5);
        }
    }
    \foreach \i in {0, 4, ..., 15} {
        \foreach \j in {0, 4, ..., 15} {
            \foreach \k in {0, 1, 2, 3} {
                \draw [thick, gray] (\i + \k, 15 - \j) -- (\i + \k, 12 - \j);
            }
            \foreach \k in {0, 1, 2} {
                \draw [thick, gray] (\i + \k, 12 - \j) -- (\i + \k + 1, 15 - \j);
            }
        }
    }
    \foreach \i in {0, 4, ..., 15} {
        \foreach \j in {0, 8} {
            \draw [thick, gray] (\i + 3, 12 - \j) -- (\i, 12 - \j - 1);
        }
    }
    
    \foreach \i in {0, 8} {
        \foreach \j in {0, 8} {
            \draw [thick, gray] (\i + 3, 8 - \j) -- (\i + 4, 15 - \j);
        }
    }
    \foreach \i in {0, 8} {
        \draw [thick, gray] (\i + 7, 8) -- (\i, 7);
    }
\end{tikzpicture}
    \caption{{\em Tensorization of a $16 \times 16$ matrix $K(i, j)$ into a six-dimensional tensor via a truncated binary tree construction (the same one as in Figure~\ref{fig:tensorization}). Simultaneously going up the row and column index trees interlaces the tensor indices, creating a tensor $\widetilde{\mathcal{K}}^{\{3\}}(\overline{i_1, i_2}, \overline{j_1, j_2}, i_3, j_3, i_4, j_4)$ of size $4^{\times 2} \times 2^{\times 4}$. The linear order through this tensor is more spatially local than the standard column major order of the original matrix. We can improve the quality of this spatial locality by increasing the depth of the truncated binary tree, which leads to better compression ratios by leveraging local smoothness of the underlying function.}}
    \label{fig:interlacing}
\end{figure}
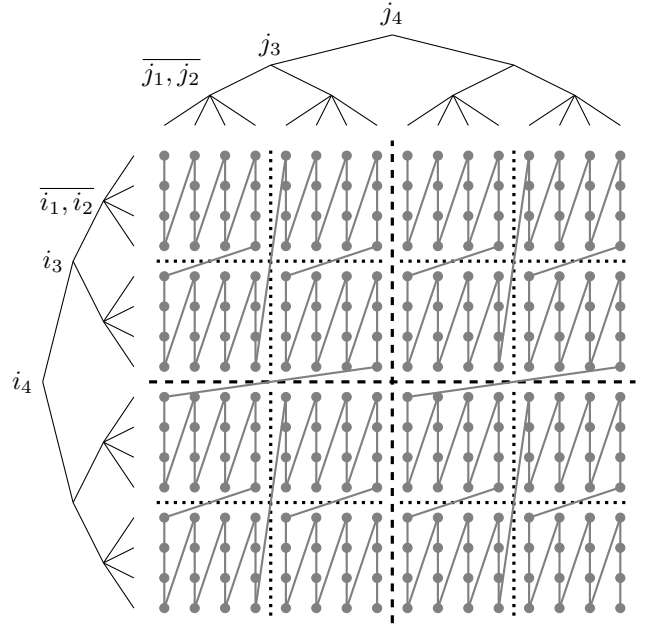

We now consider the bivariate function
\begin{equation}
    \label{eq:bivariate}
    k_\delta(x, y) = \log\frac{1}{|x - y| + \delta}, \quad \delta > 0
\end{equation}
which becomes sharper as $x \to y$ and $\delta \to 0$. We sample this function uniformly in the domain $(0, 1)^2$ for a fixed value of $\delta$ to construct the matrix $K_\delta$ of size $2^d \times 2^d$ and entries
\begin{equation}
    K_\delta(i, j) = k_\delta\left(\frac{i - 1/2}{2^d}, \frac{j - 1/2}{2^d}\right), \quad 1 \leq i, j \leq 2^d.
\end{equation}
We can tensorize this matrix, as in the univariate case, by constructing one binary tree per dimension (see Figure~\ref{fig:interlacing}). This leads to the following $2l$-dimensional tensor:
\begin{equation}
    \begin{split}
        \mathcal{K}_\delta^{\{l\}}(&\overline{i_1, \ldots, i_{d - l + 1}}, i_{d - l + 2}, \ldots, i_d, \\
        &\overline{j_1, \ldots, j_{d - l + 1}}, j_{d - l + 2}, \ldots, j_d) = \\
        & \qquad K_\delta(\overline{i_1, \ldots, i_d}, \overline{j_1, \ldots, j_d}).
    \end{split}
\end{equation}
If we try to compress this tensor directly, the resulting compression ratios will be very poor. However, if we interlace the indices of the tensor (going up row and column binary trees simultaneously),
\begin{equation}
    \begin{split}
        \widetilde{\mathcal{K}}_\delta^{\{l\}}(&\overline{i_1, \ldots, i_{d - l + 1}}, \overline{j_1, \ldots, j_{d - l + 1}}, \\
        & i_{d - l + 2}, j_{d - l + 2}, \ldots, i_d, j_d) = \\
        & \qquad K_\delta(\overline{i_1, \ldots, i_d}, \overline{j_1, \ldots, j_d}),
    \end{split}
\end{equation}
then the resulting tensor is highly compressible.

\begin{table}[t]
    \centering
    \begin{tabular}{c r@{} @{}c@{} @{}l r@{} @{}c@{} @{}l r@{} @{}c@{} @{}l}
        \toprule
                           & \multicolumn{9}{c}{Tensorization Level} \\
        \cmidrule(lr){2-10}
        $\tau_\text{relFrob}$          & \multicolumn{3}{c}{$l = 6$} & \multicolumn{3}{c}{$l = 7$} & \multicolumn{3}{c}{$l = 8$} \\
        \midrule
                           & \multicolumn{9}{c}{Compression Ratio} \\
        \cmidrule(lr){2-10}
        $10^{-2}$   & $1.6$ & $\times$ & $10^2$   & $4.9$ & $\times$ & $10^2$   & $1.0$ & $\times$ & $10^3$   \\
        $10^{-5}$   & $1.0$ & $\times$ & $10^2$   & $2.9$ & $\times$ & $10^2$   & $4.6$ & $\times$ & $10^2$   \\
        $10^{-8}$   & $7.5$ & $\times$ & $10^1$   & $2.0$ & $\times$ & $10^2$   & $3.1$ & $\times$ & $10^2$   \\
        $10^{-11}$  & $6.3$ & $\times$ & $10^1$   & $1.6$ & $\times$ & $10^2$   & $2.2$ & $\times$ & $10^2$   \\
        $10^{-14}$  & $5.4$ & $\times$ & $10^1$   & $1.3$ & $\times$ & $10^2$   & $1.7$ & $\times$ & $10^2$   \\
        \midrule
                           & \multicolumn{9}{c}{Reconstruction Error} \\
        \cmidrule(lr){2-10}
        $10^{-2}$   & $7.5$ & $\times$ & $10^{-1}$   & $7.7$ & $\times$ & $10^{-1}$   & $9.6$ & $\times$ & $10^{-1}$  \\
        $10^{-5}$   & $1.1$ & $\times$ & $10^{-3}$   & $1.3$ & $\times$ & $10^{-3}$   & $1.1$ & $\times$ & $10^{-3}$  \\
        $10^{-8}$   & $9.8$ & $\times$ & $10^{-7}$   & $1.0$ & $\times$ & $10^{-6}$   & $1.0$ & $\times$ & $10^{-6}$  \\
        $10^{-11}$  & $1.5$ & $\times$ & $10^{-9}$   & $1.5$ & $\times$ & $10^{-9}$   & $1.1$ & $\times$ & $10^{-9}$  \\
        $10^{-14}$  & $5.9$ & $\times$ & $10^{-12}$  & $6.0$ & $\times$ & $10^{-12}$  & $7.4$ & $\times$ & $10^{-12}$ \\
        \bottomrule
    \end{tabular}
    \caption{{\em Compression ratios and reconstruction errors in constructing TT compressed representation $\widehat{K}_\delta$ of the $2^d \times 2^d$ sampled kernel matrix $K_\delta$ with $d = 10$ and $\delta = 10^{-5}$ at various levels $l$ and TT-SVD tolerance $\tau_\text{relFrob}$. We use $\epsilon = \| K_\delta - \widehat{K}_\delta \|_2$ as a measure of reconstruction error. We note that the compression ratio increases as we increase the level of tensorization $l$ and decreases as we reduce TT-SVD relative tolerance and ask for a more accurate reconstruction.}}
    \label{tab:integral_compression}
\end{table}

Tensorized indices $i_k$ and $j_k$ correspond to bisecting each of the square subdomains horizontally and vertically, respectively. Arranging these indices together means we alternate between horizontal (row) and vertical (column) bisection; if we merged $i_k$ and $j_k$, this produces a quad-tree which is standard for 2D domain decomposition.

Alternatively, we may also observe that the column-major (Fortran) order imposed on the entries of permuted tensor $\widetilde{\mathcal{K}}_\delta^{\{l\}}$ obtained after interlacing the indices is more \emph{spatially coherent} than the column-major order of $\mathcal{K}_\delta^{\{l\}}$ before the interlacing --- points that are close by in the linear order of $\widetilde{\mathcal{K}}_\delta^{\{l\}}$ are also close by in the original two-dimensional domain (see Figure~\ref{fig:interlacing}). We note that as we increase $l$, the level of spatial coherence also increases. This allows higher correlation between the columns of the corresponding unfolding matrices, leading to better performance for TT compression.

In Table~\ref{tab:integral_compression} we demonstrate the compression ratios and reconstruction errors when compressing this tensorized kernel matrix at various tensorization levels $l$. We note that as we increase tensorization level up to $l = d$, the compression ratio increases.

\subsubsection{Compression of DEM Particle Data}

The previous examples pertain to compressing structured datasets obtained by sampling univariate and multi-variate functions in Cartesian domains. In this context, we demonstrated the equivalence between tensorizing dimensions of the sampled dataset and a form of hierarchical domain decomposition. By increasing the level of tensorization, we leveraged correlations in the data across multiple levels of the hierarchy, resulting in increased compression ratios.

We made two key observations. First, that ordering of the data set and of its dimensions after tensorization in a spatially coherent way is crucial; otherwise, the compression algorithm will perform poorly. Second, if a function has localized discontinuities or sharp features, we should not expect compression ratios as high as those for smooth data.

The main challenge to overcome in adapting this approach to particle level DEM simulation data is that this data is less structured --- particles can move around to potentially any location in three-dimensional space as the simulation progresses. Thus, if we initially impose a spatially coherent order for our dataset (typically including position, velocity, or force data for each of the particles over all timesteps), this coherence can degrade drastically over the course of the simulation, sometimes even in a small number of timesteps.

In the context of granular media simulations, there is a large range of phenomena for which we may assume that particles that are close by in the physical three-dimensional space will share relatively similar velocities and experience somewhat similar forces. Otherwise the movement and deformation of the granular media can change rapidly and the material is likely to fracture or fall apart.

\begin{figure}[t]
    \centering
    \includegraphics[width=0.75\linewidth]{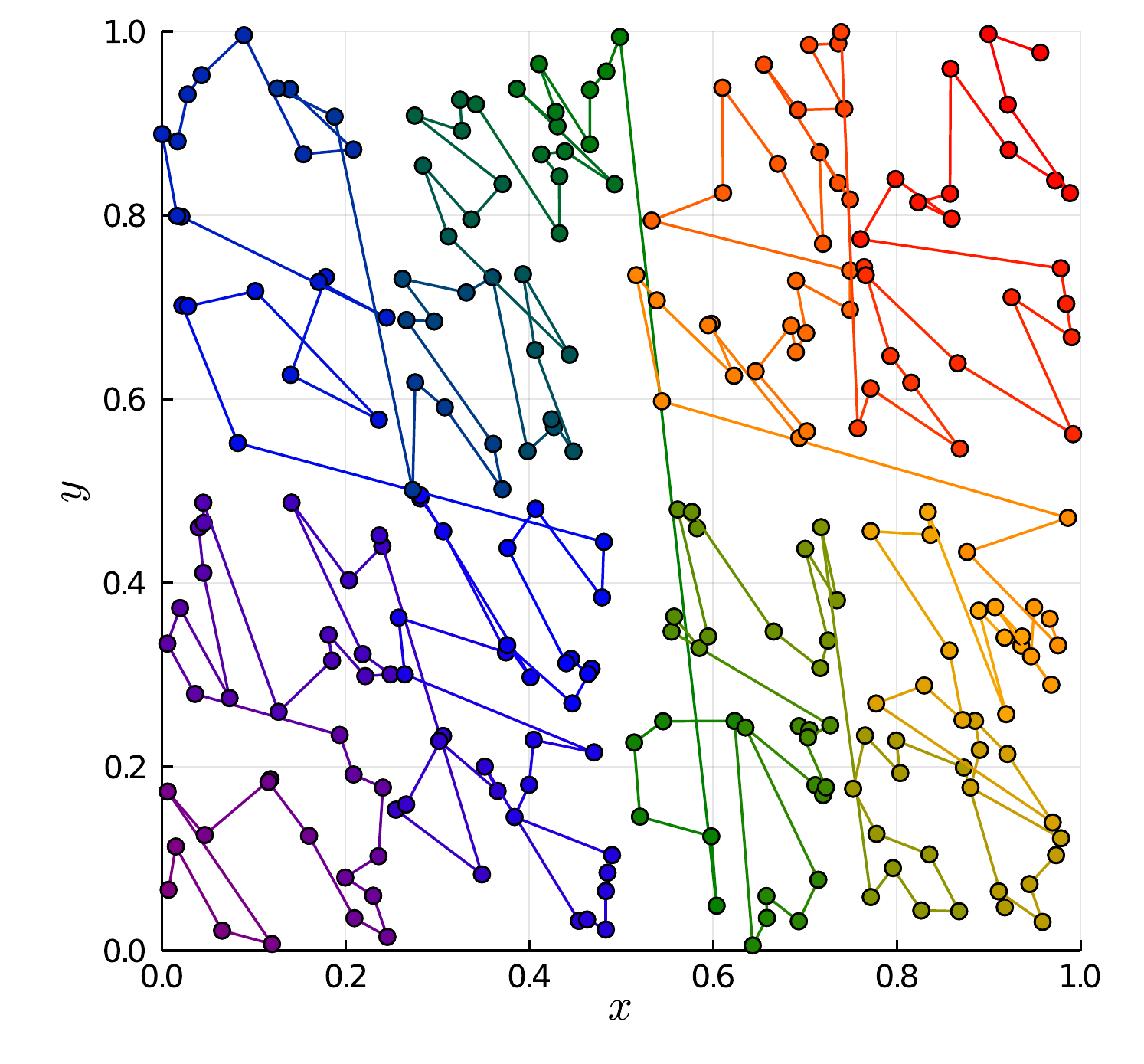}
    \caption{{\em Demonstration of spatial coherence in Morton ordering in two-dimensional Euclidean space. 256 points are generated randomly and the bit-string length $b$ was chosen base on minimal distance between the points. The color indicates the sorting using Morton indices; purple to red indicates increasing indices.}}
    \label{fig:morton-id}
\end{figure}

Since spatial coherence can be assumed to degrade slowly, we may for a given time interval $[t_0,t_f]$ exploit some of this correlation by imposing spatially coherent linear order through the physical space by utilizing particle position data at the first timestep $t_0$, see Figure~\ref{fig:morton-id}. We partition our datasets into regular time intervals and compress it using our QTT approach; we generally observe that tuning interval size provides moderate improvements for the overall compression ratio.

There are several ways we can impose a spatially coherent linear order to high-dimensional Cartesian spaces. In our numerical experiments, we use the Morton Z-order space-filling curve because of its ease of implementation. Given any DEM simulation dataset, we can always shift and rescale the simulation domain to ensure all the particles reside inside $[0, 1]^3$. Then, expanding each particle's center coordinates $(x, y, z)$ using $b$-length binary expansion
\begin{align}
    x & \approx (0.x_1 \cdots x_b)_2, \quad x_k \in \{0, 1\} \\
    y & \approx (0.y_1 \cdots y_b)_2, \quad y_k \in \{0, 1\} \\
    z & \approx (0.z_1 \cdots z_b)_2, \quad z_k \in \{0, 1\}
\end{align}
and interlacing these bit-strings we compute the Morton index for the coordinate
\begin{equation}
    \text{MortonId}_b(x, y, z) = (x_1 y_1 z_1 \cdots x_b y_b z_b)_2 \enspace .
\end{equation}
The bit-string length $b$ is chosen such that $2^{-b}$ is smaller than the extent of any of the DEM particles (e.g.\ the radius of the smallest spherical particle); this ensures the Morton index for the particle is unique.

Once we reorder the particles of a given dataset, we end up with a three-dimensional tensor of size $n_t \times n_p \times n_c$, where $n_t$ is the number of timesteps, $n_p$ is the number of particles and $n_c$ is the number of components in the variable under consideration (e.g.\ $n_c = 3$ for linear velocity). We then tensorize the time dimension and particle dimension separately and compress the resulting high-dimensional datasets.

We note that the number of particles, timesteps or variable components are usually not powers of 2 in practical applications. While this was a central theme in our development of the tensorization idea in the structured data examples, it is entirely unnecessary. As long as we can factorize these dimension sizes as product of small primes, we can easily construct high-dimensional tensors. Finally, we can also copy over data from the one of the particles (usually, the last particle in the Morton order) or the last timesteps a few times to pad the sizes of the data tensors.

\subsection{Streaming Data Compression}

One of the disadvantages of the direct TT-SVD procedure presented in Algorithm~\ref{alg:tt-svd} is the SVD operation itself. Given a $m \times n$ matrix with $m \geq n$, computing its full SVD decomposition (including the singular vectors) costs $O(m n^2)$ floating point operations. Thus, for a tensor with high number of degrees of freedom (e.g. high dimensionality), performing the SVD itself becomes a challenge.

To overcome this high computational cost, we take advantage of the properties of the TT decomposition; in particular, we utilize the fact that arithmetic operations are efficient in TT format. This allows us to break the data tensor into smaller pieces, compress them separately, and finally combine them to create the full compressed representation. This reduces the problem size for the SVD algorithm, and thereby, the overall computational cost of the data compression. Additionally, this approach can be easily adapted for compressing streaming data as long as we know the full tensor sizes in advance.

\subsubsection{Concatenation along Existing Dimension}
\label{sec:streaming_old_dim}

Many simulations have a natural time-like parameter, and the dataset grows along this dimension as the simulation progresses. In this setup, it is natural to split the data tensor along the streaming dimension, compress each piece individually, and combine the reduced representations to obtain the final compressed dataset. This incremental approach was developed in \cite{liu2018incremental, liu2020scalable}; we briefly review it here. In the next section, we adapt this approach for tensorized data.

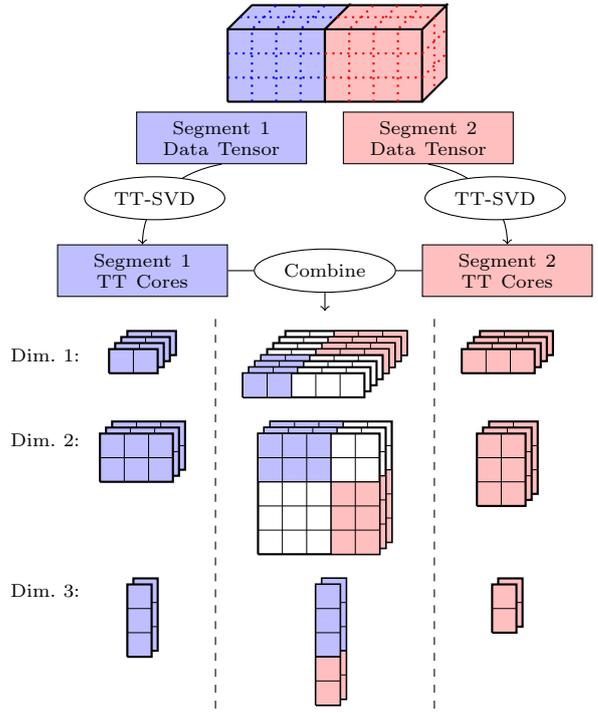
\begin{figure}[t]
    \centering

\begin{tikzpicture}[scale=0.16]
    \begin{scope}[shift={(0, 9)}]
        \begin{scope}[shift={(-8, 0)}]
            \draw[thick,fill=blue!25] (0, 0) -- (8, 0) -- (8, 6) -- (0, 6) -- (0, 0);
            \draw[thick,fill=blue!25] (8, 0) -- (10, 2) -- (10, 8) -- (8, 6) -- (8, 0);
            \draw[thick,fill=blue!25] (0, 6) -- (8, 6) -- (10, 8) -- (2, 8) -- (0, 6);
            \foreach \x in {2, 4, 6} {
                \draw[thick,blue,dotted] (\x, 0) -- (\x, 6) -- (\x + 2, 8);
            }
            \foreach \y in {2, 4} {
                \draw[thick,blue,dotted] (0, \y) -- (8, \y) -- (10, \y + 2);
            }
            \foreach \z in {1} {
                \draw[thick,blue,dotted] (\z, \z + 6) -- (\z + 8, \z + 6) -- (\z + 8, \z);
            }
        \end{scope}
        \begin{scope}[shift={(0, 0)}]
            \draw[thick,fill=red!25] (0, 0) -- (8, 0) -- (8, 6) -- (0, 6) -- (0, 0);
            \draw[thick,fill=red!25] (8, 0) -- (10, 2) -- (10, 8) -- (8, 6) -- (8, 0);
            \draw[thick,fill=red!25] (0, 6) -- (8, 6) -- (10, 8) -- (2, 8) -- (0, 6);
            \foreach \x in {2, 4, 6} {
                \draw[thick,red,dotted] (\x, 0) -- (\x, 6) -- (\x + 2, 8);
            }
            \foreach \y in {2, 4} {
                \draw[thick,red,dotted] (0, \y) -- (8, \y) -- (10, \y + 2);
            }
            \foreach \z in {1} {
                \draw[thick,red,dotted] (\z, \z + 6) -- (\z + 8, \z + 6) -- (\z + 8, \z);
            }
        \end{scope}
    \end{scope}
    \begin{scope}[shift={(0, -12)}]
        \node[draw, text width=2cm, align=center, fill=blue!25] (S1) at (-8.5, 18) {\scriptsize Segment 1 \\[-0.5em] Data Tensor};
        \node[draw, text width=2cm, align=center, fill=red!25] (S2) at (8.5, 18) {\scriptsize Segment 2 \\[-0.5em] Data Tensor};

        \node[draw, text width=2cm, align=center, fill=blue!25] (C1) at (-15, 7) {\scriptsize Segment 1 \\[-0.5em] TT Cores};
        \node[draw, text width=2cm, align=center, fill=red!25] (C2) at (15, 7) {\scriptsize Segment 2 \\[-0.5em]  TT Cores};

        \draw[->] (S1.south) to[out=190, in=90] (C1.north);
        \draw[->] (S2.south) to[out=350, in=90] (C2.north);

        \node[draw, ellipse] (C) at (0, 7) {\scriptsize Combine};
        \node[draw, ellipse, fill=white] at (-14, 13) {\scriptsize TT-SVD};
        \node[draw, ellipse, fill=white] at (14, 13) {\scriptsize TT-SVD};

        \draw (C1.east) to (C.west);
        \draw (C2.west) to (C.east);

        \draw[->] (C.south) to ++(0, -1.5);

        \draw[dashed] (-9, 3) -- (-9, -30);
        \draw[dashed] (9, 3) -- (9, -30);

        \node at (-23, 0) {\scriptsize Dim.\ 1:};
        \node at (-23, -7) {\scriptsize Dim.\ 2:};
        \node at (-23, -19.5) {\scriptsize Dim.\ 3:};

        \begin{scope}[shift={(-16.25, 0)}]
            \foreach \z in {0, 1, 2, 3} {
                \begin{scope}[shift={(-0.5 * \z, -0.5 * \z)}]
                    \draw[thick,fill=blue!25] (0, 0) -- (4, 0) -- (4, 2) -- (0, 2) -- (0, 0);
                    \foreach \x in {2} {
                        \draw[very thin] (\x, 0) -- (\x, 2);
                    }
                \end{scope}
            }
        \end{scope}
        \begin{scope}[shift={(12.75, 0)}]
            \foreach \z in {0, 1, 2, 3} {
                \begin{scope}[shift={(-0.5 * \z, -0.5 * \z)}]
                    \draw[thick,fill=red!25] (0, 0) -- (6, 0) -- (6, 2) -- (0, 2) -- (0, 0);
                    \foreach \x in {2, 4} {
                        \draw[very thin] (\x, 0) -- (\x, 2);
                    }
                \end{scope}
            }
        \end{scope}
        \begin{scope}[shift={(-3.25, 0)}]
            \foreach \z in {0, 1, 2, 3} {
                \begin{scope}[shift={(-0.5 * \z, -0.5 * \z)}]
                    \draw[thick,fill=white] (0, 0) -- (10, 0) -- (10, 2) -- (0, 2) -- (0, 0);
                    \draw[fill=red!25] (4, 0) -- (10, 0) -- (10, 2) -- (4, 2) -- (4, 0);
                    \foreach \x in {2, 4, 6, 8} {
                        \draw[very thin] (\x, 0) -- (\x, 2);
                    }
                \end{scope}
            }
            \foreach \z in {4, 5, 6, 7} {
                \begin{scope}[shift={(-0.5 * \z, -0.5 * \z)}]
                    \draw[thick,fill=white] (0, 0) -- (10, 0) -- (10, 2) -- (0, 2) -- (0, 0);
                    \draw[fill=blue!25] (0, 0) -- (4, 0) -- (4, 2) -- (0, 2) -- (0, 0);
                    \foreach \x in {2, 4, 6, 8} {
                        \draw[very thin] (\x, 0) -- (\x, 2);
                    }
                \end{scope}
            }
        \end{scope}
    \end{scope}
    \begin{scope}[shift={(0, -27.5)}]
        \begin{scope}[shift={(-17.5, 6)}]
            \foreach \z in {0, 1, 2} {
                \begin{scope}[shift={(-0.5 * \z, -0.5 * \z)}]
                    \draw[thick,fill=blue!25] (0, 0) -- (6, 0) -- (6, 4) -- (0, 4) -- (0, 0);
                    \foreach \x in {2, 4} {
                        \draw[very thin] (\x, 0) -- (\x, 4);
                    }
                    \foreach \y in {2} {
                        \draw[very thin] (0, \y) -- (6, \y);
                    }
                \end{scope}
            }
        \end{scope}
        \begin{scope}[shift={(13.5, 4)}]
            \foreach \z in {0, 1, 2} {
                \begin{scope}[shift={(-0.5 * \z, -0.5 * \z)}]
                    \draw[thick,fill=red!25] (0, 0) -- (4, 0) -- (4, 6) -- (0, 6) -- (0, 0);
                    \foreach \x in {2} {
                        \draw[very thin] (\x, 0) -- (\x, 6);
                    }
                    \foreach \y in {2, 4} {
                        \draw[very thin] (0, \y) -- (4, \y);
                    }
                \end{scope}
            }
        \end{scope}
        \begin{scope}[shift={(-4.5, 0)}]
            \foreach \z in {0, 1, 2} {
                \begin{scope}[shift={(-0.5 * \z, -0.5 * \z)}]
                    \draw[thick,fill=white] (0, 0) -- (10, 0) -- (10, 10) -- (0, 10) -- (0, 0);
                    \draw[fill=blue!25] (0, 6) -- (6, 6) -- (6, 10) -- (0, 10) -- (0, 6);
                    \draw[fill=red!25] (6, 0) -- (10, 0) -- (10, 6) -- (6, 6) -- (6, 0);
                    \foreach \x in {2, 4, 6, 8} {
                        \draw[very thin] (\x, 0) -- (\x, 10);
                    }
                    \foreach \y in {2, 4, 6, 8} {
                        \draw[very thin] (0, \y) -- (10, \y);
                    }
                \end{scope}
            }
        \end{scope}
    \end{scope}
    \begin{scope}[shift={(0, -40.5)}]
        \begin{scope}[shift={(-15.75, 4)}]
            \foreach \z in {0, 1} {
                \begin{scope}[shift={(-0.5 * \z, -0.5 * \z)}]
                    \draw[thick,fill=blue!25] (0, 0) -- (2, 0) -- (2, 6) -- (0, 6) -- (0, 0);
                    \foreach \y in {2, 4} {
                        \draw[very thin] (0, \y) -- (2, \y);
                    }
                \end{scope}
            }
        \end{scope}
        \begin{scope}[shift={(14.25, 6)}]
            \foreach \z in {0, 1} {
                \begin{scope}[shift={(-0.5 * \z, -0.5 * \z)}]
                    \draw[thick,fill=red!25] (0, 0) -- (2, 0) -- (2, 4) -- (0, 4) -- (0, 0);
                    \foreach \y in {2} {
                        \draw[very thin] (0, \y) -- (2, \y);
                    }
                \end{scope}
            }
        \end{scope}
        \begin{scope}[shift={(-0.25, 0)}]
            \foreach \z in {0, 1} {
                \begin{scope}[shift={(-0.5 * \z, -0.5 * \z)}]
                    \draw[thick,fill=white] (0, 0) -- (2, 0) -- (2, 10) -- (0, 10) -- (0, 0);
                    \draw[fill=blue!25] (0, 4) -- (2, 4) -- (2, 10) -- (0, 10) -- (0, 4);
                    \draw[fill=red!25] (0, 0) -- (2, 0) -- (2, 4) -- (0, 4) -- (0, 0);
                    \foreach \y in {2, 4, 6, 8} {
                        \draw[very thin] (0, \y) -- (2, \y);
                    }
                \end{scope}
            }
        \end{scope}
    \end{scope}
\end{tikzpicture}
    \caption{{\em Schematic of streaming TT compression of a $8 \times 3 \times 2$ tensor dataset using two equal segments along the first dimension. We first construct the TT representations of the two segments of the data (lower left and lower right of the figure) and combine them (lower middle) using appropriate zero padding (indicated by white TT core entries). Finally we run TT-rounding to optimize the TT ranks.}}
    \label{fig:streaming}
\end{figure}

In Figure~\ref{fig:streaming} we present a schematic of the incremental TT decomposition. We first partition the data tensor $\mathcal{X}$ along the first dimension into smaller segments $\mathcal{X}^{(1)}$ and $\mathcal{X}^{(2)}$ and construct their TT representations:
\begin{align}
    \mathcal{X}^{(1)}(i_1, \ldots, i_d) = X_1^{(1)}(i_1) \cdots X_d^{(1)}(i_d), \\
    \mathcal{X}^{(2)}(i_1, \ldots, i_d) = X_1^{(2)}(i_1) \cdots X_d^{(2)}(i_d).
\end{align}
Then for each dimension $1 \leq k \leq d$ we construct associated TT core slices $X_k(i_k)$ from the segment TT cores $X_k^{(j)}(:)$ for $j = 1, 2$ with appropriate zero padding; for more details on this process, see \cite{liu2018incremental, liu2020scalable}. Finally, this concatenation process usually inflates the TT ranks, so we use TT-rounding to optimize the compression ratio of the full tensor.

This scheme can be easily generalized to handle arbitrary number of data segments. Further, beyond compressing streaming data, it can also be used to compress large datasets in parallel: each segment can be compressed independently of the others, and using a binary tree structure can improve the efficiency of the recombination step even further.

\subsubsection{Stacking along New Dimension}
\label{sec:streaming_new_dim}

The approach for growing a dataset along an existing dimension proposed in \cite{liu2018incremental, liu2020scalable} is ill-suited for situations where the physical dimension (e.g.\ time) is tensorized. In this case, it is more advantageous to create new dimensions in order to preserve the hierarchical structure of the tensorized indices. We achieve this with a simple modification to the streaming compression approach as follows. Consider a tensor $\mathcal{X} \in \mathbb{R}^{n_1 \times \cdots \times n_d \times l}$ where $l$ is a small number. 
\begin{enumerate}
    \item We first compress the $d$-dimensional tensors
    \begin{equation}
        \mathcal{X}^{(k)}(i_1, \ldots, i_d) = \mathcal{X}(i_1, \ldots, i_d, k)
    \end{equation}
    individually for $1 \leq k \leq l$ using the TT-SVD algorithm with tolerance $\tau_\text{relFrob}$. Denote the compressed tensor as $\mathcal{Y}^{(k)}$ with TT cores $\mathcal{Y}_1^{(k)}, \ldots, \mathcal{Y}_d^{(k)}$.

    \item Next, define the TT cores as:
    \begin{align}
        \widehat{Y}_1(i_1) &=
        \begin{bmatrix}
            Y_1^{(1)}(i_1) & \cdots & Y_1^{(l)}(i_1)
        \end{bmatrix} \\
        \widehat{Y}_s(i_s) &=
        \begin{bmatrix}
            Y_s^{(1)}(i_s) & & \\
            & \ddots & \\
            & & Y_s^{(l)}(i_s)
        \end{bmatrix}, \; 1 < s \leq d \\
        \widehat{Y}_{d + 1}(i_{d + 1}) &= e_{i_{d + 1}}
    \end{align}
    where $e_i$ is the $i$-th column of the $l \times l$ identity matrix. Then by direct computation, we have
    \begin{equation}
        \widehat{\mathcal{Y}}(i_1, \ldots, i_{d + 1}) \approx \mathcal{X}(i_1, \ldots, i_{d + 1})
    \end{equation}
    where $\widehat{\mathcal{Y}}$ is the TT tensor with the newly constructed cores.

    \item Finally, we run the cores through the TT-rounding algorithm with tolerance $\tau_\text{round}$ and create the optimized TT tensor $\mathcal{Y}$.
\end{enumerate}

As before, this process can also be used to compress large datasets in parallel.

\subsubsection{Error Analysis}

The error analyses of both approaches to streaming data compression are exactly the same. For ease of cross-referencing the algorithm steps with the sources of errors, we focus on the ``dataset stacking'' approach given in Section~\ref{sec:streaming_new_dim}: at step~1, the compressed tensors $\mathcal{Y}^{(k)}$ satisfy
\begin{equation}
    \| \mathcal{X}^{(k)} - \mathcal{Y}^{(k)} \|_F \leq \tau_\text{relFrob} \| \mathcal{X}^{(k)} \|_F \enspace .
\end{equation}
Now, after stacking the tensor segments in step~2, the overall Frobenius error is given by
\begin{equation}
    \begin{split}
        \| \mathcal{X} - \widehat{\mathcal{Y}} \|_F^2
        &= \sum_{k = 1}^l \| \mathcal{X}^{(k)} - \mathcal{Y}^{(k)} \|_F^2 \\
        &\leq \tau_\text{relFrob}^2 \sum_{k = 1}^l \| \mathcal{X}^{(k)} \|_F^2 \\
        &= \tau_\text{relFrob}^2 \| \mathcal{X} \|_F^2
    \end{split}
\end{equation}
and it follows by the triangle inequality
\begin{equation}
    \| \widehat{\mathcal{Y}} \|_F \leq \| \mathcal{X} \|_F + \| \mathcal{X} - \widehat{\mathcal{Y}} \|_F \leq (1 + \tau_\text{relFrob}) \| \mathcal{X} \|_F
\end{equation}
Finally, after TT rounding in step~3, we have
\begin{equation}
    \| \widehat{\mathcal{Y}} - \mathcal{Y} \|_F \leq \tau_\text{round} \| \widehat{\mathcal{Y}} \|_F
\end{equation}
and once again applying the triangle inequality, we obtain
\begin{equation}
    \begin{split}
        & \quad\; \| \mathcal{X} - \mathcal{Y} \|_F \\
        &\leq \| \mathcal{X} - \widehat{\mathcal{Y}} \|_F + \| \widehat{\mathcal{Y}} - \mathcal{Y} \|_F \\
        &\leq \tau_\text{relFrob} \| \mathcal{X} \|_F + \tau_\text{round} \| \widehat{\mathcal{Y}} \|_F \\
        &\leq \tau_\text{relFrob} \| \mathcal{X} \|_F + \tau_\text{round} (1 + \tau_\text{relFrob}) \| \mathcal{X} \|_F \\
        &= \tau_\text{relFrob}^\text{new} \| \mathcal{X} \|_F \enspace ,
    \end{split}
\end{equation}
with $\tau_\text{relFrob}^\text{new} = \tau_\text{relFrob} + \tau_\text{round} + \tau_\text{relFrob} \tau_\text{round}$. Clearly, by adjusting the tolerance parameters $\tau_\text{relFrob}$ and $\tau_\text{round}$, we can achieve any target overall relative Frobenius error tolerance $\tau_\text{relFrob}^\text{new}$.

\section{Numerical Results}
\label{sec:results}

We test the performance of our tensor-train data compression approach on two kinds of data: derived DEM datasets consisting of stress-strain data obtained at reference points as part of hybrid finite element and discrete element (FE-DE) simulations \cite{yamashita2020parallelized}, and raw DEM datasets consisting of position and velocity data of the particles in the granular media simulations \cite{de2019scalable, wasfy2019understanding}. In all the experiments, we compute the following metrics:
\begin{itemize}
    \item \textbf{Compression Ratio:} For a $n_1 \times \cdots \times n_d$ tensor can be represented with TT ranks $\{r_0, \ldots, r_d\}$
    \begin{equation}
        \text{compression ratio} = \frac{n_1 \cdots n_d}{r_0 n_1 r_1 + \cdots + r_{d - 1} n_d r_d}
    \end{equation}
    is a measure of the quality of the compression.

    \item \textbf{Reconstruction Error:} For a $d$-dimensional tensor $\mathcal{X}$ and its reconstruction $\mathcal{Y}$, we use the normalized root mean squared error
    \begin{equation}
        \label{eq:nrmse}
        \epsilon_{\text{nRMSE}} = \frac{\sqrt{\frac{1}{n_{\mathcal{X}}} \sum_{i} (\mathcal{X}(i) - \mathcal{Y}(i))^2}}{x_\text{max} - x_\text{min}}
    \end{equation}
    as a measure of error incurred due to the compression. Here the sum ranges over all entries of the tensor, $x_{\text{max}}$ and $x_{\text{min}}$ are the maximum and minimum entries, and $n_{\mathcal{X}}$ is the total number of entries.
\end{itemize}

We use the open-source python library \texttt{ttpy}\footnote{\texttt{https://github.com/oseledets/ttpy}} for computing the TT decomposition of data tensors and recompressing inflated TT cores. All of the experiments presented in this section were carried out on a Dell XPS laptop with 2.5~GHz Intel Core i9 and 32~GB of RAM.

\subsection{Stress-Strain Data from Triaxial Test Simulation}
\label{sec:triaxial}

\begin{figure}[t]
    \centering
    
    \includegraphics[width=0.5\linewidth]{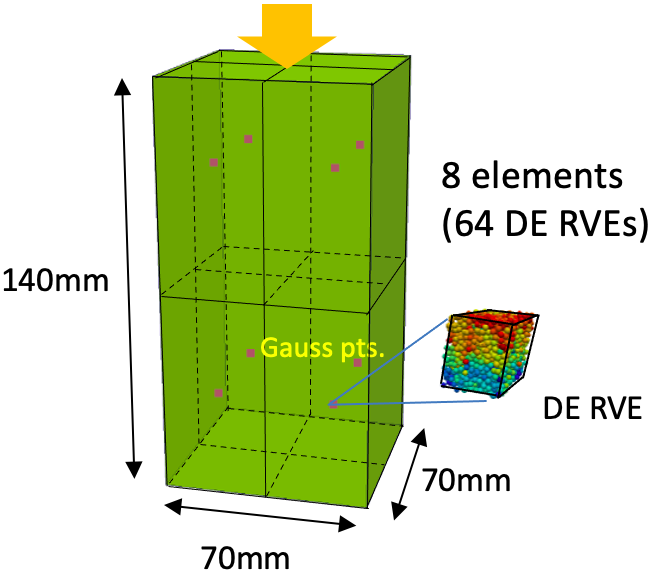}
    \caption{{\em Setup for stress-strain dataset from coupled finite element/discrete element (FE-DE) simulation of the triaxial soil compression simulation. A load is applied to the top of a column of soil confined by a constant pressure.}}
    \label{fig:triaxial-setup}

    \includegraphics*[width=\linewidth]{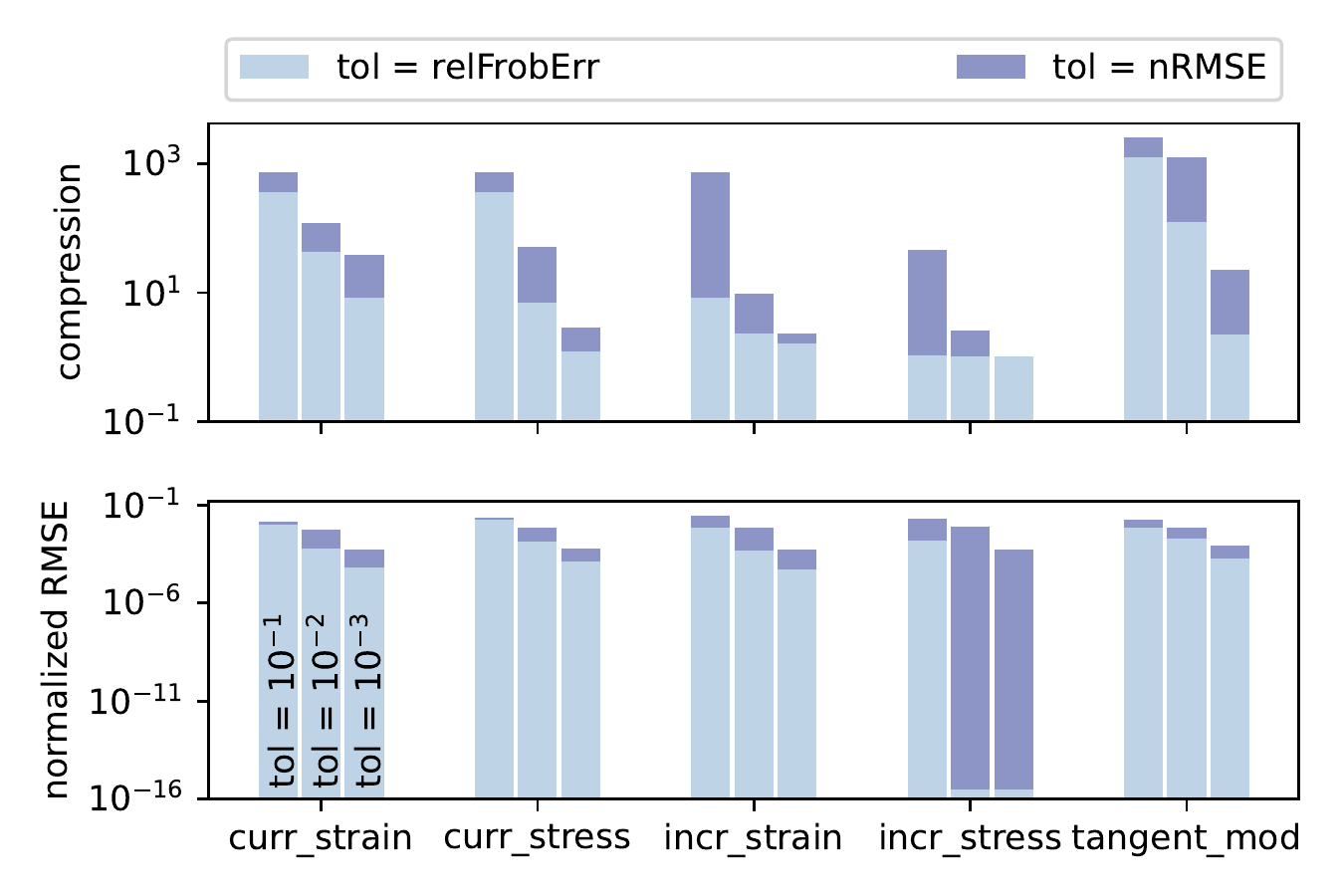}
    \caption{{\em Compression ratios and reconstruction errors of different variables in the triaxial compression test dataset for various reconstruction tolerances $\{10^{-1}, 10^{-2}, 10^{-3}\}$. The light purple bar plots correspond to using relative Frobenius error as the reconstruction metric; we see from the bottom panel that this approach usually overestimates the rank in order to meet the relative Frobenius error criteria built into the TT-SVD algorithm. Targeting the more forgiving normalized RMSE criteria directly leads to better compression ratios while still respecting the error guarantees, as indicated by the dark purple bar plots.}}
    \label{fig:triaxial-compress}
\end{figure}

In this simulation, a column of soil is confined by a constant pressure and a load is applied to the top of the column. The column is modeled using a $2 \times 2 \times 2$ finite element model (FEM) grid. Using the numerical algorithms developed in~\cite{yamashita2020parallelized}, the stress-strain response of each element is estimated from reference volume element (RVE) computations at 8 Gaussian points within the element. Each RVE computation consists of a small-scale DEM simulation of $\sim 1000$ particles (see Figure~\ref{fig:triaxial-setup}).

For each of the 64 RVEs and for two confining pressures, the simulation generates current strain, current stress, incremental strain, incremental stress, and tangent moduli data---with $n_c = 6, 6, 6, 6, 21$ components, respectively---for 2028 timesteps; this leads to five tensors of size $2 \times 2028 \times 64 \times n_c$.

We compress each of these tensors using TT-SVD with specified relative tolerances $\tau_\text{relFrob} \in \{10^{-1}, 10^{-2}, 10^{-3}\}$ and plot the compression ratios and reconstruction errors in Figure~\ref{fig:triaxial-compress} using light purple bar plots. We note that as we demand a more accurate representation by decreasing the relative tolerance parameter of the TT-SVD algorithm, we achieve more accurate representations at the cost of lower compression ratios.

We also note that in some cases the compressed dataset obtained from TT-SVD has inflated ranks, leading to low compression ratios. This is due to the the fact that the relative Frobenius error criteria built into the algorithm which is more conservative than the normalized RMSE reconstruction error metric. This is particularly evident for the incremental stress variable, where TT-SVD tolerance of $\tau_\text{relFrob} = 10^{-2}$ leads to a normalized RMSE of $\sim 10^{-15}$.

To remedy this, we target the normalized RMSE reconstruction error criteria directly when specifying the tolerance of the TT-SVD algorithm: the normalized RMSE \eqref{eq:nrmse} can be rewritten as
\begin{equation}
    \epsilon_{\text{nRMSE}} = \frac{1}{x_{\text{max}} - x_{\text{min}}} \frac{\| \mathcal{X} \|_F}{\sqrt{n_{\mathcal{X}}}} \underbrace{\frac{\| \mathcal{X} - \mathcal{Y} \|_F}{\| \mathcal{X} \|_F}}_{\epsilon_\text{relFrob}}.
\end{equation}
Given this direct connection between the normalized RMSE and relative Frobenius error, we can specify a target normalized RMSE tolerance $\tau_\text{nRMSE}$ and derive an equivalent relative Frobenious error tolerance \begin{equation}
    \tau_\text{relFrob} = (x_\text{max} - x_\text{min}) \frac{\sqrt{n_\mathcal{X}}}{\| \mathcal{X} \|_F} \tau_\text{nRMSE}
\end{equation}
to use in TT-SVD. This can also be used for streaming data compression since Frobenius norms are cheap to compute in TT format \cite{oseledets2011tensor} and we can easily keep track of the minimum/maximum entries of each data segment during the recombination steps.

We compress the variables in our triaxial dataset using this new target error criteria, and plot the corresponding compression ratios/reconstruction errors in Figure~\ref{fig:triaxial-compress} using the dark purple bar plots --- we observe drastic improvements in the compression ratios for some of the variables.

\begin{figure}[t]
    \centering
    \includegraphics[width=\linewidth]{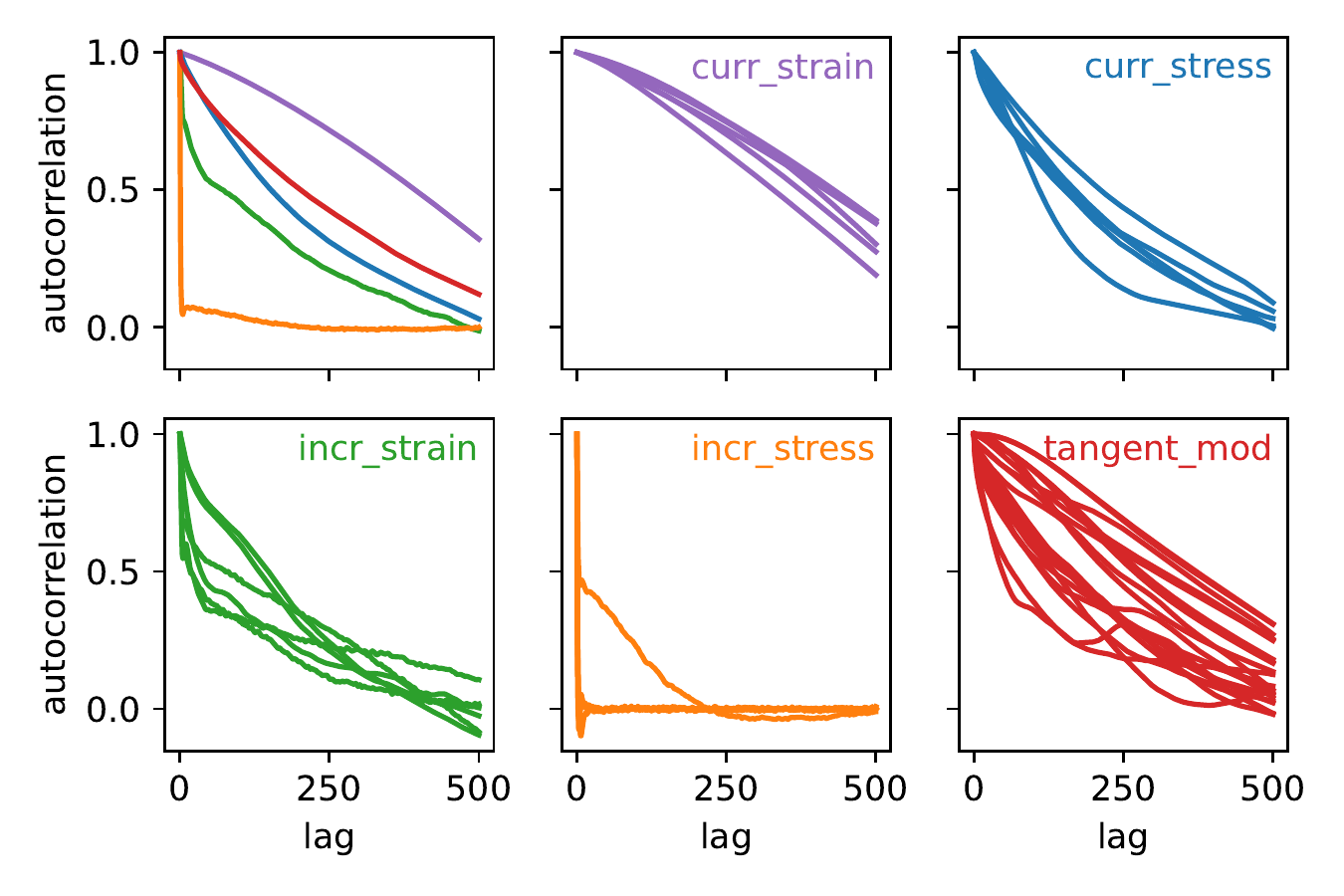}
    \caption{{\em Autocorrelations of the triaxial compression simulation datasets along the time dimension and averaged over the confining pressure and RVE dimensions. The top-left panel plots the autocorrelations averaged over the variable components, while the other panels plot the autocorrelations for each individual component. The rate of decay of the autocorrelation with increasing lag is a measure of self-similarity/redundancy of the dataset --- a slow decay corresponds to highly correlated variables. This leads to better performance of the compression algorithm for the current strain and tangent moduli variables, and very poor performance for the incremental stress variable.}}
    \label{fig:triaxial-autocorr}
\end{figure}

Finally, we observe that some variables are more compressible than the others. Intuitively, if a dataset is highly uncorrelated (or noisy), then compressing it is harder. This can be corroborated by analyzing the autocorrelations; given a timeseries $\{y_1, \ldots, y_n\}$, the $k$-lag autocorrelation is defined as
\begin{equation}
    \text{autocorrelation}(k) = \frac{\frac{1}{n - k} \sum_{i = 1}^{n - k} (y_i - \bar{y}) (y_{i + k} - \bar{y})}{\frac{1}{n} \sum_{i = 1}^n (y_i - \bar{y})^2}
\end{equation}
with
\begin{equation}
    \bar{y} = \frac{1}{n} \sum_{i = 1}^n y_i \enspace .
\end{equation}
Then, by definition
\begin{equation}
    \text{autocorrelation}(k) =
    \begin{cases}
        1 & \text{if} \quad k = 0, \\
        0 & \text{if} \quad k \geq n.
    \end{cases}
\end{equation}
and the rate of decay for $0 < k < n$ is an indicator of self-similarity in the timeseries. In Figure~\ref{fig:triaxial-autocorr} we plot the autocorrelations of the variables computed along the time dimension, and averaged over the pressure and RVE dimensions. We note that,
\begin{itemize}
    \item The autocorrelations for the current stress and tangent moduli variables decay relatively slowly. This is an indicator that the datasets are self-similar, corresponding to higher compression ratios.

    \item The autocorrelation for the incremental stress variable decays very quickly. Without any redundancy to exploit, the compression algorithm performs poorly on this variable on the dataset.
\end{itemize}

\subsection{Stress-Strain Data from Soil-Wheel Simulation}
\label{sec:soilwheel}

\begin{figure}[t]
    \centering
    \includegraphics[width=0.5\linewidth]{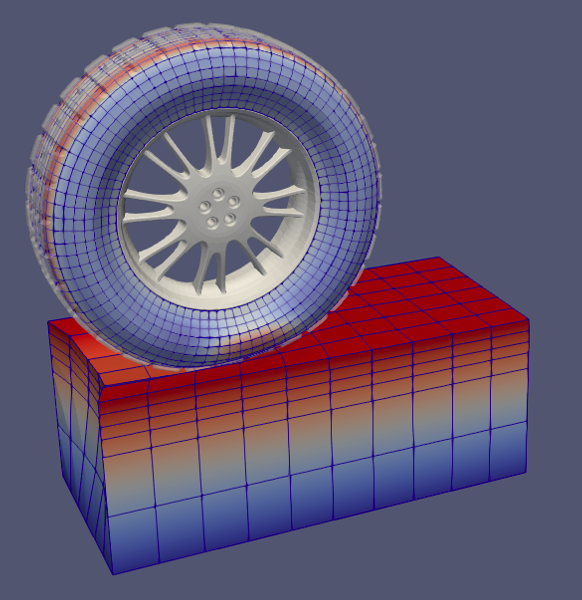}
    \caption{{\em Setup for the stress-strain datasets from the soil-wheel FE-DE simulation \cite{yamashita2020parallelized}. A deformable wheel is rolled on top of a soil patch at constant velocity, and the soil-patch is modeled with a $10 \times 6 \times 6$ grid.  Within each element, 8 RVEs are placed on a $2 \times 2 \times 2$ quadrature grid, and stress strain data is generated from all of these RVEs.}}
    \label{fig:soilwheel-setup}
\end{figure}

\begin{table}[h!]
    \centering

    \begin{tabular}{c c c c}
        \toprule
                      & \multicolumn{3}{c}{Tolerance $\tau_\text{nRMSE}$} \\
        \cmidrule(lr){2-4}
        Dataset       & $10^{-1}$     & $10^{-2}$     & $10^{-3}$     \\
        \midrule
                      & \multicolumn{3}{c}{Compression Ratio}                              \\
        \cmidrule(lr){2-4}
        Curr.\ Strain & $6.3 \times 10^4$    & $8.0 \times 10^2$    & $8.2 \times 10^1$    \\
        Curr.\ Stress & $4.0 \times 10^4$    & $7.3 \times 10^2$    & $8.1 \times 10^0$    \\
        Incr.\ Strain & $2.0 \times 10^4$    & $1.2 \times 10^3$    & $8.9 \times 10^0$    \\
        Incr.\ Stress & $3.0 \times 10^4$    & $1.9 \times 10^3$    & $1.4 \times 10^0$    \\
        Tangent Mod.  & $4.7 \times 10^5$    & $7.0 \times 10^4$    & $1.8 \times 10^1$    \\
        \cmidrule(lr){2-4}
        Overall       & $7.3 \times 10^4$    & $2.5 \times 10^3$    & $7.6 \times 10^0$    \\
        \midrule
                      & \multicolumn{3}{c}{Normalized RMSE}                               \\
        \cmidrule(lr){2-4}
        Curr.\ Strain & $6.2 \times 10^{-2}$ & $8.2 \times 10^{-3}$ & $6.0 \times 10^{-4}$ \\
        Curr.\ Stress & $2.5 \times 10^{-2}$ & $8.4 \times 10^{-3}$ & $6.1 \times 10^{-4}$ \\
        Incr.\ Strain & $2.0 \times 10^{-2}$ & $8.6 \times 10^{-3}$ & $6.1 \times 10^{-4}$ \\
        Incr.\ Stress & $1.5 \times 10^{-2}$ & $8.2 \times 10^{-3}$ & $6.1 \times 10^{-4}$ \\
        Tangent Mod.  & $2.1 \times 10^{-2}$ & $8.0 \times 10^{-3}$ & $9.1 \times 10^{-4}$ \\
        \bottomrule
    \end{tabular}
    \caption{{\em Compression ratios and reconstruction errors from TT compression of soil-wheel simulation datasets. We compress each dataset with three different relative tolerances $\tau_\text{nRMSE} \in \{10^{-1}$, $10^{-2}, 10^{-3}\}$ and record the compression ratio and the reconstruction errors.}}
    \label{tbl:soilwheel-compression}
\end{table}

In this simulation, a wheel is rolled on top of a soil patch at constant velocity. Both the wheel and the soil is modeled using FEM methods; in particular, the soil patch is discretized using a $10 \times 6 \times 6$ grid. The response of each of the elements is computed from DE simulations at $2 \times 2 \times 2$ Gaussian quadrature grids --- this leads to a total of $2880$ reference volume elements (RVEs). Stress-strain data --- consisting of current strain, current stress, incremental strain, incremental stress and tangent moduli with $n_c = 6$, 6, 6, 6, and 36 components respectively --- is generated for 602 timesteps and 5 different patches. This creates 9-dimensional tensors of size $5 \times 602 \times 10 \times 6 \times 6 \times 2 \times 2 \times 2 \times n_c$.

\begin{figure*}
    \centering
    \includegraphics[width=\linewidth]{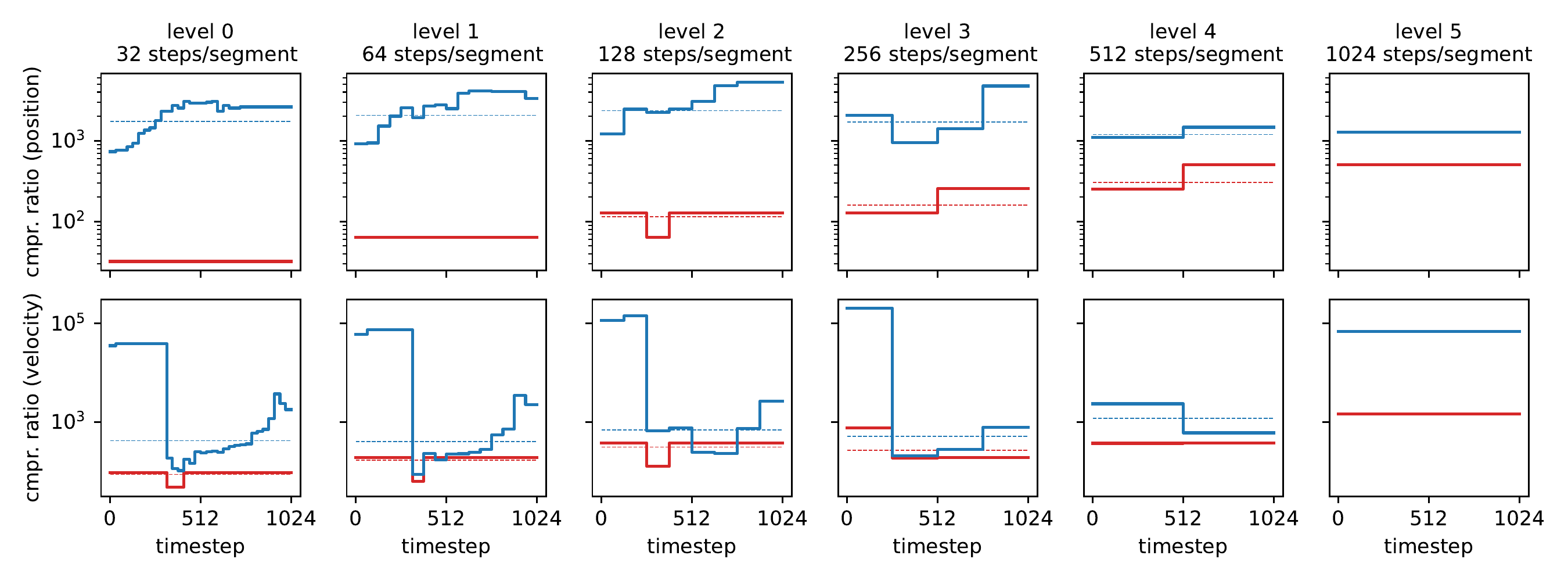}
    \caption{{\em Streaming compression of the position and velocity variables from the sedimentation dataset with maximum normalized RMSE of $10^{-1}$. \textbf{Leftmost Column:} We compress the position and velocity variables in 32 segments (each segment containing data from 32 timesteps of the simulation) and plot the compression ratios for the original 3D data segments (red) and tensorized 16D data segments (blue). We also add the overall compression ratio (averaged over all segments) using the dashed lines. \textbf{Subsequent Columns:} We take the TT tensors from the previous level and combine them, two at a time, to construct compressed representations of the data over longer time segments (e.g.\ the second level TT tensors have 64 timesteps per segment, ..., the final level TT tensors have a single segment with all 1024 timesteps). As with the top level figures, we plot the compression ratios for 3D datasets (red), tensorized datasets (blue) and overall compression ratios (dashed). We observe that compression ratios from the tensorized datasets are almost always significantly higher than those from the natural 3D datasets.}}
    \label{fig:sediment-str}
\end{figure*}

In Table~\ref{tbl:soilwheel-compression}, we summarize the results from compressing each of these tensors at different values of the normalized RMSE tolerances $\tau_\text{nRMSE} \in \{10^{-1}$, $10^{-2}, 10^{-3}\}$. We observe high overall compression ratios ($\sim 10^4$) with reconstruction errors smaller than $10^{-1}$. This reduces the size of the dataset from 3.9~GB to mere 102~KB.

Finally, we note that these compression ratios are an order of magnitude better than those corresponding to the triaxial test simulation dataset. We attribute this to the larger size of the soil-wheel dataset --- we can exploit the redundancies of a bigger dataset more effectively and obtain larger compression ratios.

\subsection{Particle Data from Sediment Simulation}
\label{sec:sediment}

In this simulation, $2^{14}$ spherical particles are released from rest inside a rectangular box, and the particles fall under the influence of the gravity. During the course of the simulation, the particles collide with each other and the inner walls of the box, until they settle down. The simulation is performed with a timestep size of $10^{-3}$~s for $2^{10}$ timesteps. This creates 3-dimensional datasets for positions and velocities of the particles of size $1024 \times 16384 \times 3$. We compress this tensor 32 timesteps at a time --- i.e.\ we compress 32 segments, each with size $32 \times 16384 \times 3$, for either variables.

We also use the tensorization approach to increase the dimensionality of the variable segments --- we first sort the particles using the Morton order computed from the first snapshot of each segment, then reshape the resulting 3-dimensional tensors to 20-dimensional tensors of size $2^{\times 19} \times 3$, and finally compress these high-dimensional tensor segments. In the leftmost column of Figure~\ref{fig:sediment-str} we compare the compression ratios obtained from these two approaches for the position and velocity variables.

We note that in general, the tensorized dataset segments achieve significantly higher compressibility than their natural dataset segment counterparts. The exception to this is for the velocity variable between timesteps 300 and 800 approximately. We remark that during this period, the particles start to collide with each other and the interior walls after their initial free-fall phase. Once the particles settle down, the compression ratio starts to increase again.

\begin{table}
    \centering

    \begin{tabular}{c c c c}
        \toprule
        Level & \#steps/segment & Position & Velocity \\
        \midrule
        0 &   32 & $1.7 \times 10^3$ & $4.2 \times 10^2$ \\
        1 &   64 & $2.1 \times 10^3$ & $4.0 \times 10^2$ \\
        2 &  128 & $2.4 \times 10^3$ & $6.9 \times 10^2$ \\
        3 &  256 & $1.7 \times 10^3$ & $5.2 \times 10^2$ \\
        4 &  512 & $1.2 \times 10^3$ & $1.2 \times 10^3$ \\
        5 & 1024 & $1.3 \times 10^3$ & $6.8 \times 10^4$ \\
        \bottomrule
    \end{tabular}
    \caption{{\em Overall compression ratios of the tensorized variables from the sedimentation dataset with normalized RMSE smaller than $10^{-1}$ at various levels of the streaming compression. The second column records the number of timesteps per chunk at each level --- for instance, at level 2 there are 8 chunks of the datasets are compressed separately, each containing 128 timesteps per chunk. These overall compressions correspond to the blue dashed lines in Figure~\ref{fig:sediment-str}.}}
\end{table}

\begin{figure*}
    \centering
    \includegraphics[width=0.87\linewidth]{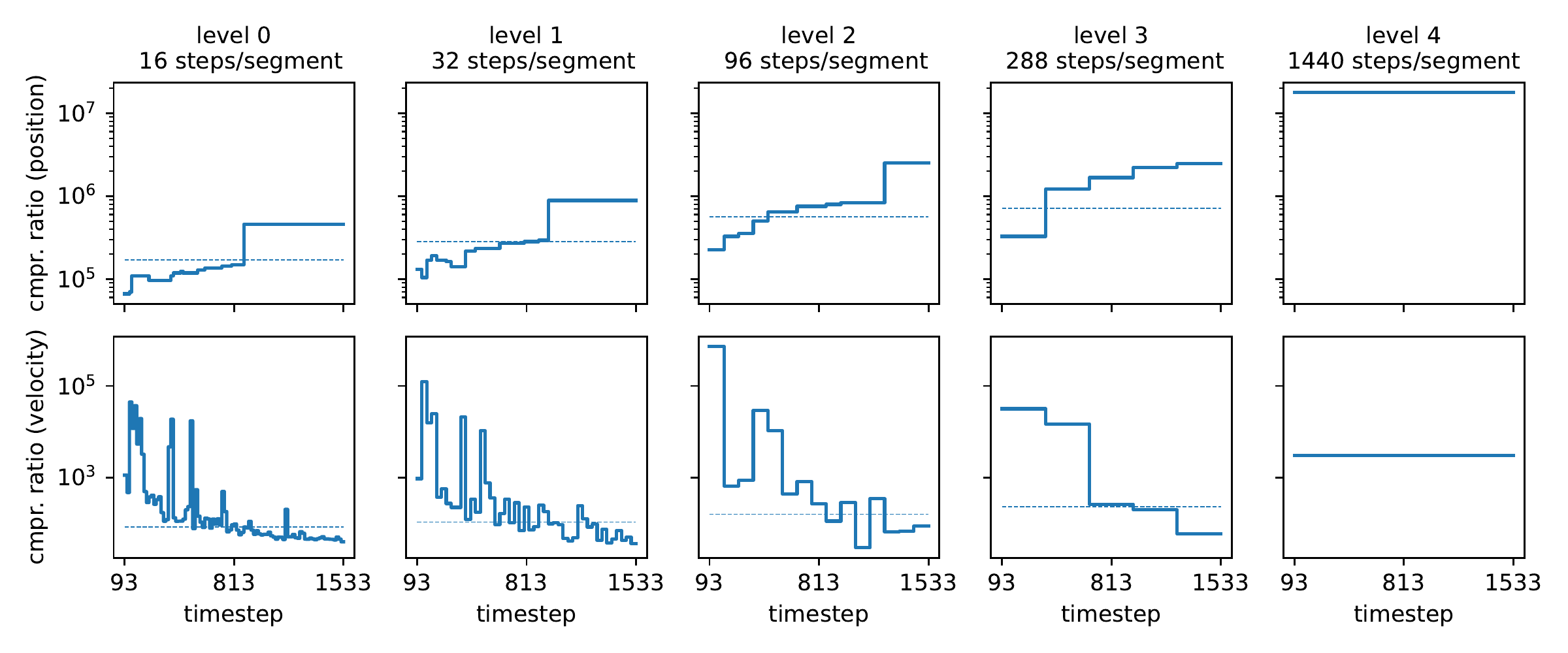}
    \caption{{\em Streaming compression of the position and velocity variables from the soil-vehicle simulation dataset with maximum normalized RMSE of $10^{-1}$. \textbf{Leftmost Column:} We compress the position and velocity variables in 90 segments and record the compression ratio. Each segment contains tensorized data from 16 timesteps of the simulation (we start compressing the data at timestep 93 --- this corresponds to the second snapshot in Figure~\ref{fig:soilvehicle-setup}). We also add the overall compression raitos using the dashed lines. \textbf{Subsequent Columns:} We take the TT tensors from the previous levels and combine them to construct compressed representation of the data over progressively bigger time segments, until we collect all timesteps into single TT tensor per variable.}}
    \label{fig:soilvehicle-str}
\end{figure*}

Next, we apply the streaming compression setup to the compressed chunks of this sediment simulation dataset. We combine the TT compressed datasets, preserving the tensorization structure, as outlined in Section~\ref{sec:streaming_new_dim}. In Figure~\ref{fig:sediment-str} we plot the compression ratios of these incrementally constructed compressed datasets. We observed the same compressibility features as with the individual data segments: an initially high compression ratio, followed by a drop when the particles collide and the data is nonsmooth, and finally an increase in the compression ratio when the particles settle down.

\subsection{Particle Data from Soil-Vehicle Simulation}
\label{sec:soilvehicle}

In this simulation, a vehicle is driven on a soil bed consisting of 488376 spherical particles. During the course of the simulation, the patch moves with the vehicle --- particles at the tail end are discarded, and an equal number of particles are generated at the front of the soil patch; a roller is used to settle these generated particles (see Figure~\ref{fig:soilvehicle-setup}). Position and velocity data are saved for 1533 timesteps at 0.03~s intervals. We discard the data from the initial 93 timesteps --- during this time initial soil bed is constructed via sedimentation and the vehicle is placed on top of it. Thus we obtain $1440 \times 488376 \times 3$ sized datasets for the position and velocity.

During the simulation, the soil patch moves with the vehicle --- this is achieved by deleting particles that are far behind the vehicle, and generating an equal number of soil particles a fixed distance ahead of the vehicle. A roller is then used to smooth over the front of the soil patch. The velocity data of the particles in this region are artificially high due to this setup, so we zero out the velocities of all particles in the first 2~m of the soil patch for all timesteps.

To tensorize these modified datasets, we first sort the particles using the Morton ordering at each timestep, then copy the last particle in this sorted order over and over again until we have a total of 500000 particles. This allows us to factorize the number of particles in terms of small prime numbers. We then compress the dataset 16 timesteps at a time. After tensorizing the time dimension as well, the resulting 16-dimensional tensorized datasets have size $2^{\times 9} \times 5^{\times 6} \times 3$. In the leftmost column of Figure~\ref{fig:soilvehicle-str} we plot the compression ratios from the TT compression of these position and velocity datasets with normalized RMSE smaller than $10^{-1}$.

\begin{table}
    \centering

    \begin{tabular}{c c c c}
        \toprule
        Level & \#steps/segment & Position & Velocity \\
        \midrule
        0 &   16 & $1.7 \times 10^5$ & $8.5 \times 10^1$ \\
        1 &   32 & $2.8 \times 10^5$ & $1.1 \times 10^2$ \\
        2 &   96 & $5.6 \times 10^5$ & $1.6 \times 10^2$ \\
        3 &  288 & $7.2 \times 10^5$ & $2.4 \times 10^2$ \\
        4 & 1440 & $1.8 \times 10^7$ & $3.1 \times 10^3$ \\
        \bottomrule
    \end{tabular}
    \caption{{\em Overall compression ratios of the tensorized variables from the soil-vehicle dataset with maximum normalized RMSE of $10^{-1}$ at various levels of the streaming compression. The second column records the number of timesteps per chunk at each level --- for instance, at level 2 there are 15 chunks of the datasets are compressed separately, each containing 96 timesteps per chunk. These overall compressions correspond to the blue dashed lines in Figure~\ref{fig:soilvehicle-str}.}}
    \label{tbl:soilvehicle}
\end{table}

Next, we combine these data segments using our streaming compression approach, while preserving the hierarchical structure of the tensorized datasets, until all the timesteps are collected in a single compressed tensor. We show the compression ratios of the dataset segments during this process in Figure~\ref{fig:soilvehicle-str}. We also record the overall compression ratios at each level in Table~\ref{tbl:soilvehicle}. We note that the position variable is highly compressible --- with final compression ratio of $1.8 \times 10^7$ with maximum normalized RMSE of $10^{-1}$. This reduces the original 15.7~GB dataset (after discarding approximately 1~GB data from initial 93 timesteps) to less than 1~KB.

While the velocity variable is not as compressible, we nonetheless reduce it from 15.7~GB to approximately 5.2~MB with maximum normalized RMSE of $10^{-1}$. This relatively poor performance of the TT compression is likely due to the non-smooth nature of the velocity data --- the large stepsize ($\Delta t = 0.03$~s) between saving the simulation snapshots degrades the correlation between the velocity data far quicker than that of the position data. This in turn leads to low compressibility of the saved velocity data.

\section{Conclusions}

In this paper, we proposed a tensor-train (TT) method for scientific data compression, tailoring our approach to datasets defined in unstructured spatial locations, as is typical of those generated by DEM simulations. For this purpose, we borrowed key concepts from applications of quantized TT to fast numerical solvers for partial differential equations (PDEs):
\begin{itemize}
    \item \emph{Spatial Locality}. DEM variables correspond to particle or collision locations; their values are likely to be similar for nearby locations in 3D space (especially in dense granular media). We used Morton ordering to impose a linear ordering that improves spatial locality in these datasets.
    
    \item \emph{Hierarchical Structure}. Tensorization of a dataset along space and time dimensions imposes a hierarchical data decomposition. Tensor-train can exploit low rank structure at different levels of this hierarchy. 
\end{itemize}

Throughout our numerical experiments, we consistently found that DEM datasets are highly compressible in TT format, and further, that the tensorization of the Morton-ordered datasets greatly improves their compressibility. 

We implemented a streaming compression scheme that incrementally compresses and merges segments of these datasets preserving data tensorization. This was instrumental to compress large DEM datasets with the robust TT-SVD algorithm while avoiding its application to the entire data set, which would have been prohibitively expensive.

Since TT data compression is lossy by design, we explored two criteria for characterizing reconstruction errors in our experiments. Relative Frobenius error is the standard criteria built into the TT-SVD algorithm, and it compares the RMSE against the average size of the tensor entries. Normalized RMSE on the other hand compares the RMSE against the range of the data tensor values. We generally found normalized RMSE to be a more forgiving criteria compared to relative Frobenius error. Hence, larger compression ratios can be achieved for this target, especially for nonsmooth variables.

Utilizing the locality and hierarchical structure in the datasets, and using normalized RMSE as reconstruction error criteria, we were able to achieve very high compression ratios for both the raw DEM position-velocity data and derived stress-strain data. For instance, we were able to compress the soil-wheel stress-strain dataset (Section~\ref{sec:soilwheel}) by a factor of $7.3 \times 10^4$, and the soil-vehicle position and velocity datasets (Section~\ref{sec:soilvehicle}) by factors of $1.8 \times 10^7$ and $3.1 \times 10^3$, respectively. 

In general, we observe that our approach allows us to greatly compress large datasets by exploiting redundancies in the data, reducing their size from GBs to a few MBs or even KBs. Once in compressed format, we can take advantage of one of the main features of tensor factorizations like the TT: data access, visualization and analysis can all be performed efficiently in TT format. This enables efficient data manipulation on machines with limited computing power and storage capacity once the simulation datasets have been compressed.

%

\section*{Acknowledgments}

We thank Hiroyuki Sugiyama, Tamer Wasfy and Hiroki Yamashita for several helpful discussions pertaining to this work and sharing their DEM simulation datasets.

We acknowledge support from the Automotive Research Center at the University of Michigan (UM) in accordance with Cooperative Agreement W56HZV-19-2-0001 with U.S. Army DEVCOM Ground Vehicle Systems Center. This research was supported in part through computational resources and services provided by UM’s Advanced Research Computing.

\bibliographystyle{elsarticle-num}
\bibliography{references}

\end{document}